\documentclass[11pt,a4paper,reqno]{amsart}
\usepackage[T1]{fontenc}
\usepackage[english]{babel}
\usepackage{lmodern,amsmath,amssymb,amsthm,mathtools,mathrsfs,microtype,booktabs}
\usepackage[a4paper,margin=28mm,headheight=14pt]{geometry}
\usepackage{needspace}
\usepackage[hidelinks]{hyperref}
\hypersetup{pdftitle={Einstein contact metrics with integrable CR structure},pdfauthor={Beniamino Cappelletti-Montano}}
\numberwithin{equation}{section}
\newtheorem{theorem}{Theorem}[section]
\newtheorem{proposition}[theorem]{Proposition}
\newtheorem{lemma}[theorem]{Lemma}
\newtheorem{corollary}[theorem]{Corollary}
\theoremstyle{definition}

\theoremstyle{remark}
\newtheorem{remark}[theorem]{Remark}
\newcommand{\R}{\mathbb R}
\newcommand{\HH}{\mathcal H}
\newcommand{\EE}{\mathcal E}
\newcommand{\KK}{\mathcal K}
\newcommand{\RR}{\mathcal R}
\newcommand{\BB}{\mathsf B}
\newcommand{\LL}{\mathcal L}
\newcommand{\bB}{\mathscr B}
\newcommand{\bC}{\mathscr C}
\newcommand{\bT}{\mathscr T}
\newcommand{\bF}{\mathscr F}
\newcommand{\Id}{\operatorname{Id}}
\newcommand{\Ric}{\operatorname{Ric}}
\newcommand{\tr}{\operatorname{tr}}
\newcommand{\rank}{\operatorname{rank}}
\newcommand{\im}{\operatorname{Im}}

\newcommand{\divg}{\operatorname{div}_g}
\newcommand{\Lie}{\mathcal L}
\newcommand{\dd}{\mathrm d}

\newcommand{\cyc}{\mathop{\mathfrak S}}
\newcommand{\norm}[1]{\lVert#1\rVert}
\allowdisplaybreaks[2]
\title[Einstein contact metrics with integrable CR structure]{Einstein contact metrics with integrable CR structure}
\author{Beniamino Cappelletti-Montano}
\address{Department of Mathematics and Computer Science, University of Cagliari, Palazzo delle Scienze, Via Ospedale 72, 09124 Cagliari, Italy}
\email{b.cappellettimontano@unica.it}
\thanks{Partially supported by project ProBiki, funded by Fondazione di Sardegna.}
\date{}

\begin{document}
\begin{abstract}
We prove that every compact Einstein contact metric manifold of dimension at least
$5$ with integrable associated CR structure 
is Sasakian. We show that without the assumption of compactness this result does not hold: namely, for every odd dimension at least $5$, we  construct an
explicit noncompact Einstein contact metric manifold, with integrable CR-structure, which is not Sasakian. 
\end{abstract}
\maketitle

\section{Introduction}\label{sec:intro}

For a contact metric structure $(\phi,\xi,\eta,g)$, the symmetric
$(1,1)$-tensor field $h=\frac12\Lie_\xi\phi$ measures the failure of the Reeb
vector field to preserve the structure. Its vanishing is equivalent to the Reeb vector field being Killing, i.e. the contact metric manifold being K-contact. If the associated CR structure is integrable,
$h=0$ is also equivalent to the Sasakian condition. It is therefore
natural to ask whether the Einstein condition forces $h$ to vanish.
Boyer and Galicki proved that every compact Einstein K-contact
manifold is Sasakian (see \cite{BG}); a different proof was given by
Apostolov, Dr\u{a}ghici and Moroianu in \cite{ADM}.  Those results assume that the Reeb vector field preserves the contact metric structure $(\phi,\xi,\eta,g)$ -- i.e. the tensor field $h$ vanishes identically -- and 
 establish CR-integrability.
 In this paper we deal with 
 the complementary hypothesis: CR-integrability
is assumed, and Reeb invariance is obtained from the Einstein condition.
Namely we shall prove the following result.

\begin{theorem}\label{thm:main}
Let $(M^{2n+1},\phi,\xi,\eta,g)$ be a compact Einstein contact metric
manifold, where $n\ge2$. If its associated CR structure is integrable,
then the given contact metric structure is Sasakian.
\end{theorem}

We use the normalization $\dd\eta(X,Y)=2g(X,\phi Y)$ and formulate
integrability through the real identity
\begin{equation}\label{eq:intro-integrability}
 (\nabla_X\phi)Y=g(X+hX,Y)\xi-\eta(Y)(X+hX).
\end{equation}
The connection used in the proof is Tanno's generalized Tanaka--Webster
connection, introduced by an explicit formula in terms of $\nabla$.
Under \eqref{eq:intro-integrability} it preserves $g,\phi,\xi$ and $\eta$.
We work throughout with real tangent bundles, real traces and real
covariant derivatives.

The first step for proving Theorem \ref{thm:main}  is to establish that the eigenvalues of $h$ are constant. Namely, we prove the following general result, that does not require the assumptions of
compactness or completeness.

\begin{theorem}\label{thm:spectrum-intro}
On an Einstein contact metric manifold of dimension $2n+1\ge5$ with
integrable associated CR structure, the eigenvalues of $h$, including
their multiplicities, are constant. In particular, $h$ has constant rank
and its eigenspaces define smooth subbundles.
\end{theorem}

Since \(h\) is symmetric and anticommutes with \(\phi\), its nonzero eigenvalues occur in opposite pairs. Thus, if \(h\neq0\), its rank must be \(2r\) for some \(1\le r\le n\). We then distinguish two cases.
If \(r=n\), then \(h\) is invertible on \(\mathcal H\). In this case, working with the Tanaka--Webster connection, we define  a $1$-form $\alpha$ (see \eqref{eq:alpha}) whose differential is, up to sign, the Ricci-type $2$-form \eqref{eq:rho-def} of the Tanaka--Webster curvature, and we show that the Einstein condition forces $\alpha=-n^2\eta$. Comparing these two facts we obtain $\lambda=2n^2-2$, $\lambda$ being the Einstein constant. Since $\lambda=\Ric_g(\xi,\xi)=2n-\|h\|^2\le2n$, this is impossible when $n\ge2$, and  hence the maximal-rank case cannot occur. Notice that this argument is purely local. Whereas, if \(1\le r<n\), 
using again the Tanaka--Webster connection, we  construct a globally defined  vector field \(\mathcal V\) satisfying
\[
\operatorname{div}_g\mathcal V
=-2r(n-r)-\frac1n U,
\]
where $U$ is a non-negative function on $M$. The right-hand side is strictly negative, contradicting the divergence theorem on the compact manifold \(M\). Thus intermediate ranks are also excluded.
It follows that \(h=0\), and  the CR-integrability then  implies that the original contact metric structure is Sasakian.

Summing up, the maximal rank is excluded by a local argument, whereas compactness plays an essential role in the exclusion of the intermediate ranks. Indeed, the intermediate-rank case can occur on noncompact manifolds, even with positive Einstein constant and complete Reeb vector field. In fact we shall prove the following  result.

\begin{theorem}\label{thm:examples-intro}
For each $n\ge2$, there is a contact metric structure
$(\phi,\xi,\eta,g)$ on $\R^{2n+1}$ whose associated CR structure is
integrable and for which
\[
 \Ric_g=2(n-1)g,\qquad
 \operatorname{spec}(h)=\{1,-1,0\},
\]
with multiplicities $1,1,2n-1$, respectively. The Reeb vector field is complete,
whereas the metric is incomplete and has unbounded sectional curvature.
Moreover, no open subset of this metric is homothetic to a
Sasaki--Einstein metric.
\end{theorem}

These examples show that Theorem~\ref{thm:main} is a global rigidity
statement: neither positivity of the Einstein constant nor completeness
of the Reeb flow can replace compactness. They have the smallest possible
nonzero rank of $h$ and satisfy the negative-divergence identity above
explicitly (see the end of Section~\ref{sec:examples}). They also differ from the familiar non-Sasakian
contact metric $(\kappa,\mu)$-spaces, where $h$ has full horizontal
rank.

Finally, we point out that also the assumption on the dimension is essential in Theorem \ref{thm:main}. In fact
the flat $3$-torus recalled in
Section~\ref{sec:remarks} prevents an extension of the statement to
dimension $3$.

\bigskip

Throughout the paper, all manifolds are assumed to
be connected, smooth and without boundary.

\section{Preliminaries}\label{sec:prelim}
Let $(\phi,\xi,\eta,g)$ be a contact metric structure on a $(2n+1)$-dimensional smooth manifold $M$. This means that $\phi$ is a $(1,1)$-tensor field, $\eta$ is a contact form with Reeb vector field $\xi$, and $g$ is a Riemannian metric satisfying
\begin{gather}
 \eta(\xi)=1,\qquad \phi^2=-\Id+\eta\otimes\xi,\qquad
 g(\phi X,\phi Y)=g(X,Y)-\eta(X)\eta(Y),\label{eq:contactalg}\\
 \dd\eta(X,Y)=2g(X,\phi Y). \label{eq:contactnormalization}
\end{gather}
We refer to \cite{Blair} for more details on contact Riemannian geometry. 
Let $\HH=\ker\eta$ be the  contact distribution. 
It is well known that in any contact metric manifold $\nabla_{\xi}\xi=0$ and $\nabla_{\xi}\phi=0$, where  $\nabla$ denotes the Levi-Civita connection of $(M,g)$.

A key role in contact Riemannian geometry is played by  the $(1,1)$-tensor field $h=\frac{1}{2}{\mathcal L}_{\xi}\phi$. This endomorphism is symmetric, traceless, anti-commutes with $\phi$ and satisfies $h\xi=0$.
Moreover,  for any vector field $X$,
\begin{equation}\label{formula_h}
 \nabla_X\xi=-\phi X-\phi hX.
\end{equation}
The vanishing of $h$ is equivalent to the Reeb vector field being Killing. In this case the contact metric structure is called K-contact. A contact metric manifold is said to be Sasakian if
\begin{equation}\label{sasakian_condition}
  (\nabla_X\phi)Y=g(X,Y)\xi-\eta(Y)X,  
\end{equation}
for any $X,Y\in\Gamma(TM)$. Any Sasakian manifold is also K-contact. 

The tangent bundle of $M$ splits as the ortogonal sum $TM=\mathcal H \oplus \mathbb{R}\xi$. We will refer to vectors tangent to $\mathcal H$ as \emph{horizontal}. Putting
$J=-\phi|_{\HH}$, one has that, since $J^2=-I$ on $\HH$, the pair
$(\HH,J)$ defines an almost CR structure of hypersurface
type on $M$. We extend $J$ to $TM$ by setting $J\xi=0$,
so that $J^2=-P$, where $P=I-\eta\otimes\xi$.
Moreover, our choice of sign ensures that the Levi form
\[
L_\eta(X,Y)=\frac12 d\eta(X,JY)
\]
coincides with $g|_{\HH}$ and is therefore positive
definite.
Let $\HH^{1,0}$ denote the $i$-eigenbundle of $J$ in
$\HH\otimes\mathbb{C}$. The almost CR structure $(\HH,J)$ is
said to be \emph{integrable} if
\begin{equation}\label{CR-complex}
[\Gamma(\HH^{1,0}),\Gamma(\HH^{1,0})]
\subseteq\Gamma(\HH^{1,0}).
\end{equation}
This means that
\begin{equation}\label{CR-real}
[JX,JY]-[X,Y]=J\bigl([JX,Y]+[X,JY]\bigr)
\end{equation}
for all $X,Y\in\Gamma(\HH)$. Note that, for a contact metric structure, the vector fields $[JX,JY]-[X,Y]$ and $[JX,Y]+[X,JY]$ are necessarily horizontal, since $\dd\eta(JX,JY)=\dd\eta(X,Y)$ and $\dd\eta(JX,Y)+\dd\eta(X,JY)=0$ on $\HH$.
When \eqref{CR-complex}, or equivalently \eqref{CR-real}, holds, $(\HH,J)$ is a \emph{strictly
pseudoconvex CR structure}.

\begin{remark}
In the contact metric literature, the opposite
convention $J=\phi|_{\mathcal H}$ is also common,
with a corresponding change in the sign convention
for the Levi form. Both choices give equivalent
CR-integrability conditions.
\end{remark}

Equivalently, the CR-integrability can be expressed in terms of the covariant derivative of $\phi$ as (see \cite{Tanno}, cf.  also \cite{Blair})
\begin{equation}\label{eq:CRlevi}
 (\nabla_X\phi)Y=g(X+hX,Y)\xi-\eta(Y)(X+hX).
\end{equation}
Comparing \eqref{sasakian_condition} and \eqref{eq:CRlevi} we see that a K-contact CR-integrable manifold is necessarily Sasakian.

Let $D$ denote the generalized Tanaka--Webster
connection, introduced by Tanno in \cite{Tanno}, and defined by
\begin{equation}\label{eq:Ddefinition}
D_XY=\nabla_XY+\eta(X)\phi Y
-\eta(Y)\nabla_X\xi
+(\nabla_X\eta)(Y)\xi.
\end{equation}
Under the assumption of this paper that the almost CR structure $(\HH,J)$ is  integrable, $D$ coincides with the usual
Tanaka--Webster connection associated with the
pseudo-Hermitian structure $(\HH,J,\eta)$.
Set
\begin{equation}\label{eq:tau}
 \tau=Jh.
\end{equation}
The connection $D$ satisfies
\begin{gather}
 Dg=DJ=D\eta=D\xi=0,\label{eq:Dmetric}\\
 T^D(X,Y)=2\omega(X,Y)\xi,\qquad
 T^D(\xi,X)=\tau X,\label{eq:torsion}
\end{gather}
for any $X,Y\in\Gamma(\HH)$,
where $\omega=\tfrac12\dd\eta=g(J\,\cdot,\cdot)$.
The pseudo-Hermitian torsion $\tau$ is symmetric
and anticommutes with $J$. Throughout the paper,
we shall denote by
\[
H=\tau^2=h^2
\]
the square of the torsion endomorphism.
Thus $H$ is symmetric, nonnegative, and commutes
with $J$. Moreover, $H\xi=0$.

Using \eqref{formula_h}, the defining formula
for $D$ can equivalently be written as
\begin{equation}\label{eq:connection}
 \nabla_XY=D_XY+\eta(X)JY
 +\eta(Y)(J+\tau)X
 -g((J+\tau)X,Y)\xi.
\end{equation}
In particular, for any  $X,Y\in\Gamma(\HH)$,
\begin{equation}\label{eq:horizontalD}
 D_XY=P\nabla_XY,\qquad
 D_\xi Y=\nabla_\xi Y-JY.
\end{equation}
Let $R^g$ and $\RR=R^D$ denote the curvature endomorphisms of the connections $\nabla$ and $D$, respectively. Since $D$ preserves $g,J$
and $\xi$, each $\RR(X,Y)$ is skew-adjoint, commutes with $J$ and
annihilates $\xi$.
For horizontal arguments, the first Bianchi identity for $D$ reduces to
\begin{equation}\label{eq:bianchi-horizontal}
 \cyc_{X,Y,Z}\RR(X,Y)Z
   =2\cyc_{X,Y,Z}\omega(X,Y)\tau Z.
\end{equation}
Indeed, the covariant derivatives of $T^D$ on horizontal pairs vanish
by \eqref{eq:Dmetric}; the remaining terms are
$T^D(T^D(X,Y),Z)=2\omega(X,Y)\tau Z$.
The following algebraic decomposition will be useful.

\begin{lemma}\label{lem:curvature-decomposition}
Let $(M,\phi,\xi,\eta,g)$ be a CR-integrable contact metric manifold. 
For $U,V\in\HH$, let $(U\wedge V)Z=g(V,Z)U-g(U,Z)V$, and set
\begin{equation}\label{eq:curvature-model}
 \RR_\tau(X,Y)=JX\wedge\tau Y-JY\wedge\tau X
                      -X\wedge J\tau Y+Y\wedge J\tau X.
\end{equation}
Then $\RR_0=\RR|_{\HH}-\RR_\tau$ is an algebraic curvature tensor on
$\HH$ which commutes with $J$. In particular,
\begin{equation}\label{eq:R0-J}
 \RR_0(JX,JY)=\RR_0(X,Y),\qquad
 \RR_0(JX,Y)=-\RR_0(X,JY).
\end{equation}
If $K$ and $K_0$ are the horizontal Ricci tensors of $\RR$ and
$\RR_0$, respectively, then
\begin{equation}\label{eq:K-decomposition}
 K(X,Y)=K_0(X,Y)-2(n-1)g(J\tau X,Y),
\end{equation}
where $K_0$ is symmetric and $J$-invariant. The real $2$-form on $M$ defined by
\begin{equation}\label{eq:rho-def}
 \rho(X,Y)=\tfrac12\tr_{\HH}(J\RR(X,Y))
\end{equation}
satisfies
\begin{equation}\label{eq:rho-K0}
 \rho(X,Y)=K_0(JX,Y), 
\end{equation}
for all $X,Y\in\Gamma(\HH)$.
\end{lemma}

\begin{proof}
Every summand in \eqref{eq:curvature-model} is skew-adjoint, and
$\tau J=-J\tau$ gives $[\RR_\tau(X,Y),J]=0$.
Expansion of the cyclic sum gives
\[
 \cyc_{X,Y,Z}\RR_\tau(X,Y)Z
     =2\cyc_{X,Y,Z}\omega(X,Y)\tau Z.
\]
Thus $\RR_0$ is skew in its first two and last two arguments and satisfies
the ordinary first Bianchi identity. These properties imply that $g(\RR_0(X,Y)Z,W)
=g(\RR_0(Z,W)X,Y)$.
Commutation with $J$ then gives
\eqref{eq:R0-J} and the symmetry and $J$-invariance of $K_0$. Next, let $\{e_A\}_{A\in\{1,\ldots,2n\}}$ be a local orthonormal frame of $\HH$.  By \eqref{eq:curvature-model},
the horizontal Ricci tensor of $\RR_\tau$
is given by
\[
\begin{aligned}
\sum_{A=1}^{2n}
g(\RR_\tau(e_A,X)Y,e_A)&=\sum_A g((Je_A\wedge\tau X)Y,e_A)
-\sum_A g((JX\wedge\tau e_A)Y,e_A)\\
&\quad-\sum_A g((e_A\wedge J\tau X)Y,e_A)
+\sum_A g((X\wedge J\tau e_A)Y,e_A)\\
&=-g(J\tau X,Y)+g(J\tau X,Y)-(2n-1)g(J\tau X,Y)+g(J\tau X,Y)\\
&=-2(n-1)g(J\tau X,Y),
\end{aligned}
\]
where the four sums have been computed separately, using $g(Je_A,e_A)=0$ and $\tr_{\HH}\tau=\tr_{\HH}(J\tau)=0$. Hence
\[
K(X,Y)=K_0(X,Y)+\sum_{A=1}^{2n}g(\RR_\tau(e_A,X)Y,e_A)=K_0(X,Y)-2(n-1)g(J\tau X,Y),
\]%
which proves \eqref{eq:K-decomposition}.
Likewise,
$\tr(J\RR_\tau(X,Y))=0$. Applying the ordinary Bianchi identity to
$\RR_0(X,Y)e_A$, pairing with $Je_A$ and summing gives
\[
 \sum_A g(\RR_0(X,Y)e_A,Je_A)=-2K_0(JX,Y).
\]
Here the other two cyclic terms each equal $K_0(JX,Y)$ by
\eqref{eq:R0-J} and pair interchange. Since
$\tr(J\RR_0(X,Y))=-\sum_Ag(\RR_0(X,Y)e_A,Je_A)$,
this proves \eqref{eq:rho-K0}.
\end{proof}

Since $D$ preserves $\HH$, the endomorphism $D_X\tau$ preserves
$\HH$ for every vector field $X$ on $M$. Its horizontal contraction
defines a global horizontal vector field $\delta\tau$. In any local orthonormal
frame $\{e_A\}_{A\in\{1,\ldots,2n\}}$ of $\HH$, it is
given by
\[
\delta\tau=\sum_{A=1}^{2n}(D_{e_A}\tau)e_A.
\]
We refer to $\delta\tau$ as the horizontal
divergence of $\tau$.

\begin{lemma}\label{lem:einstein}
Let $(M,\phi,\xi,\eta,g)$ be a CR-integrable contact metric manifold. 
Suppose that $\Ric_g=\lambda g$. Then
\begin{gather}
 D_\xi\tau=-2nJ\tau,\quad D_\xi H=0,\quad
 \tr_{\HH}H=2n-\lambda,\label{eq:einstein-reeb}\\
\delta\tau=0,\label{eq:divtau}\\
 K_0=(\lambda+2)g|_{\HH},\quad
 \rho|_{\HH\times\HH}=(\lambda+2)\omega|_{\HH\times\HH}.
 \label{eq:webster-einstein}
\end{gather}
\end{lemma}

\begin{proof}
Put $B=J+\tau=\nabla\xi$. Since $B^2=-P+\tau^2$ and $\xi$ is
geodesic, its Jacobi operator $\ell X=R^g(X,\xi)\xi$ satisfies
$\ell=-\nabla_\xi B-B^2$. By \eqref{eq:horizontalD},
$\nabla_\xi B=D_\xi\tau+2J\tau$, and hence
\begin{equation}\label{eq:jacobi-tau}
 \ell=P-\tau^2-2J\tau-D_\xi\tau.
\end{equation}
Both $J\tau$ and $D_\xi\tau$ anticommute with $J$ on $\HH$ and
annihilate $\xi$, so they are traceless. Consequently,
\begin{equation}\label{eq:ricxi}
 \Ric_g(\xi,\xi)=2n-\tr_{\HH}H=2n-\norm h^2.
\end{equation}
Let $C_XY=(\nabla_X-D_X)Y$. Then the curvature difference formula is
\begin{align*}
 R^g(X,Y)Z={}&\RR(X,Y)Z+(D_XC)_YZ-(D_YC)_XZ\\
 &+C_XC_YZ-C_YC_XZ+C_{T^D(X,Y)}Z.
\end{align*}
For horizontal $X,Y,Z$, formula \eqref{eq:connection} gives
$C_XY=-g(BX,Y)\xi$, $C_X\xi=BX$ and $C_\xi Z=JZ$.
The terms involving $DC$ are vertical. Therefore
\begin{equation}\label{eq:curvature-horizontal}
 PR^g(X,Y)Z=\RR(X,Y)Z-g(BY,Z)BX+g(BX,Z)BY
                                      +2\omega(X,Y)JZ.
\end{equation}
It follows that
\begin{align*}
 \sum_Ag(R^g(e_A,X)Y,e_A)
 &=K(X,Y)+g(B^2X,Y)-2g(X,Y)\\
 &=K(X,Y)+g(HX,Y)-3g(X,Y).
\end{align*}
The remaining Ricci contribution is
$g(R^g(\xi,X)Y,\xi)=g(\ell X,Y)$. Using \eqref{eq:jacobi-tau}
and \eqref{eq:K-decomposition}, we obtain
\begin{equation}\label{eq:ric-final}
 \Ric_g(X,Y)=K_0(X,Y)-2g(X,Y)-2n g(J\tau X,Y)
                                      -g((D_\xi\tau)X,Y).
\end{equation}
The first two terms are $J$-invariant; the last two are
$J$-anti-invariant. Thus the Einstein condition implies
$K_0=(\lambda+2)g$ and $D_\xi\tau=-2nJ\tau$.
Equation \eqref{eq:rho-K0} proves the identity for $\rho$.
Differentiating $H=\tau^2$ gives $D_\xi H=0$, and
\eqref{eq:ricxi} gives its trace.
Finally, the curvature comparison formula
gives, for horizontal vector fields $X,Y$,
\[
R^g(X,Y)\xi
=(D_X\tau)Y-(D_Y\tau)X.
\]
Taking the horizontal Ricci contraction
and using $\tr_{\HH}\tau=0$, we obtain
\[
\Ric_g(X,\xi)=g(\delta\tau,X).
\]
Since $g$ is Einstein and $X$ is horizontal,
the left-hand side vanishes. Consequently,
$\delta\tau=0$, proving \eqref{eq:divtau}.
\end{proof}

\section{Spectral properties of the pseudo-Hermitian torsion}\label{sec:spectrum}

In this section we prove Theorem~\ref{thm:spectrum-intro}. Thus we assume that $(M,\phi,\xi,\eta,g)$ is a contact metric manifold with CR-integrable structure, such that $\Ric_g=\lambda g$.
 No compactness  assumption
is needed. For each integer $m\ge1$, put
\begin{equation}\label{eq:spectral-moments}
 s_m=\tr_{\HH}(H^m).
\end{equation}
These are the power sums of the eigenvalues of $H$ on the contact distribution
$\HH$. We already know, by Lemma \ref{lem:einstein} that $s_1=2n-\lambda$ is constant. Our objective is to show that all the other $s_m$ are
constant as well.

For horizontal vector fields define the tensors
\begin{align}
 t(X,Y,Z)&=g((D_Z\tau)X,Y),\label{eq:real-t}\\
 \bB(X,Y,Z)&=t(X,Y,Z)-t(JX,Y,JZ),\label{eq:real-B}\\
 \bC(X,Y,Z)&=t(X,Y,Z)+t(JX,Y,JZ),\label{eq:real-C}\\
 \bT(X,Y,Z)&=\bC(X,Y,\tau Z).\label{eq:real-T}
\end{align}
Note that, since $D_Z\tau$ is symmetric and anticommutes with $J$,
\begin{equation}\label{eq:t-symmetries}
 t(X,Y,Z)=t(Y,X,Z),\qquad t(JX,Y,Z)=t(X,JY,Z).
\end{equation}
Thus $\bB,\bC,\bT$ are symmetric in their first two arguments.
Moreover, as a straightforward consequence of \eqref{eq:real-B}--\eqref{eq:real-T}, one has
\begin{gather}
\bB(JX,Y,Z)=\bB(X,JY,Z)=\bB(X,Y,JZ),\label{eq:B-Jslots}\\
 \bC(JX,Y,Z)=\bC(X,JY,Z)=-\bC(X,Y,JZ),\label{eq:C-Jslots}\\
 \bT(JX,Y,Z)=\bT(X,JY,Z)=\bT(X,Y,JZ).\label{eq:T-Jslots}
\end{gather}

\begin{lemma}\label{lem:real-Codazzi}
For any $X,Y,Z\in\Gamma(\mathcal{H})$,
\begin{equation}\label{eq:real-mixedcurvature}
 g(\RR(\xi,X)Y,Z)=t(X,Z,Y)-t(X,Y,Z).
\end{equation}
Moreover, the tensor $\bB$ is totally symmetric.
\end{lemma}

\begin{proof}
We apply the first Bianchi identity for the
connection $D$ to the
triple $(\xi,X,Y)$, where $X,Y\in\Gamma(\HH)$.
Since $D\xi=0$, we have
$\RR(X,Y)\xi=0$, and hence
\begin{align}\label{bianchi1}
\RR(\xi,X)Y-\RR(\xi,Y)X&=(D_\xi T^D)(X,Y)
 +(D_XT^D)(Y,\xi)
 +(D_YT^D)(\xi,X)\\
&\quad+T^D(T^D(\xi,X),Y)
+T^D(T^D(X,Y),\xi)+T^D(T^D(Y,\xi),X). \nonumber
\end{align}
Using $D\omega=0$, $D\xi=0$, and the
torsion identities \eqref{eq:torsion},
we obtain
\[
\begin{aligned}
(D_\xi T^D)(X,Y)=0, \quad
(D_XT^D)(Y,\xi)=-(D_X\tau)Y,\quad
(D_YT^D)(\xi,X)=(D_Y\tau)X.
\end{aligned}
\]
Moreover, 
\[
\begin{aligned}
T^D(T^D(\xi,X),Y)
 =2\omega(\tau X,Y)\xi,\ \
T^D(T^D(X,Y),\xi)=0,\ \
T^D(T^D(Y,\xi),X)=-2\omega(\tau Y,X)\xi.
\end{aligned}
\]
Since $\tau$ is symmetric and anti-commutes with $J$, the endomorphism
$J\tau$ is symmetric. Consequently,
\[
\omega(\tau X,Y)
=g(J\tau X,Y)
=g(J\tau Y,X)
=\omega(\tau Y,X).
\]
Thus \eqref{bianchi1} becomes
\begin{equation}\label{bianchi2}
\RR(\xi,X)Y-\RR(\xi,Y)X
=(D_Y\tau)X-(D_X\tau)Y.
\end{equation}
Since $Dg=0$, the curvature endomorphism
$\RR(\xi,X)$ is skew-adjoint. Taking the
inner product with $Z$ in \eqref{bianchi2}, we obtain
\begin{equation}\label{prima}
g(\RR(\xi,X)Y,Z)
+g(\RR(\xi,Y)Z,X)=t(X,Z,Y)-t(Y,Z,X).
\end{equation}
Applying the same identity to the cyclic
permutations of $X,Y,Z$ gives
\begin{equation}\label{seconda}
g(\RR(\xi,Y)Z,X)
+g(\RR(\xi,Z)X,Y)=t(Y,X,Z)-t(Z,X,Y),
\end{equation}
and
\begin{equation}\label{terza}
g(\RR(\xi,Z)X,Y)
+g(\RR(\xi,X)Y,Z)=t(Z,Y,X)-t(X,Y,Z).
\end{equation}
Adding \eqref{prima} and \eqref{terza}, and
subtracting \eqref{seconda}, we obtain
\begin{align}\label{bianchi3}
2g(\RR(\xi,X)Y,Z)
=&t(X,Z,Y)-t(Y,Z,X)
+t(Z,Y,X)\\
&-t(X,Y,Z)
-t(Y,X,Z)+t(Z,X,Y).\nonumber
\end{align}
Since $\tau$ is symmetric and $Dg=0$,
the tensor $t$ is symmetric in its first
two arguments. Hence the right-hand side
reduces to
$2t(X,Z,Y)-2t(X,Y,Z)$.
Therefore \eqref{bianchi3} becomes
\begin{equation*}
g(\RR(\xi,X)Y,Z)
=t(X,Z,Y)-t(X,Y,Z),
\end{equation*}
as claimed. Moreover, using that  $\RR(\xi,X)$ commutes with $J$, from 
\eqref{eq:real-mixedcurvature} and  \eqref{eq:t-symmetries} it follows that
\[
 t(X,Y,Z)-t(JX,Y,JZ)
       =t(X,Z,Y)-t(JX,Z,JY).
\]
Thus $\bB(X,Y,Z)=\bB(X,Z,Y)$. 
\end{proof}

\Needspace{7\baselineskip}
\begin{lemma}\label{lem:secondcodazzi}
The tensor $\bT$ is totally symmetric.
\end{lemma}

\begin{proof}
By the definition of the curvature tensor,
the covariant derivatives of $\tau$ satisfy
\[
\begin{aligned}
&D_X(D_Y\tau)-D_Y(D_X\tau)
-D_{[X,Y]}\tau=
\RR(X,Y)\circ\tau-\tau\circ\RR(X,Y).
\end{aligned}
\]
Taking $X=\xi$ and using $D\xi=0$ and
$T^D(\xi,Z)=\tau Z$, we obtain
\[
\begin{aligned}
D_\xi(D_Z\tau)-D_{D_\xi Z}\tau
={}&D_Z(D_\xi\tau)+[\RR(\xi,Z),\tau]
-D_{\tau Z}\tau,
\end{aligned}
\]
which, using $D_\xi\tau=-2nJ\tau$ and $DJ=0$,
 becomes
\begin{equation}\label{eq:real-commutator}
\begin{aligned}
D_\xi(D_Z\tau)-D_{D_\xi Z}\tau
={}&-2nJ(D_Z\tau)+[\RR(\xi,Z),\tau]
-D_{\tau Z}\tau.
\end{aligned}
\end{equation}
Since $Dg=0$, the definition of $t$ gives
\[
\begin{aligned}
(D_\xi t)(X,Y,Z)
={}&g\bigl(
(D_\xi(D_Z\tau)-D_{D_\xi Z}\tau)X,Y
\bigr).
\end{aligned}
\]
We now take the inner product with $X,Y$
in \eqref{eq:real-commutator} and compute
its three terms separately.
First, differentiating $\tau J=-J\tau$
and using $DJ=0$, we obtain
$(D_Z\tau)J=-J(D_Z\tau)$.
Hence, by the skew-adjointness of $J$,
\[
\begin{aligned}
-2n\,g(J(D_Z\tau)X,Y)
&=2n\,g((D_Z\tau)JX,Y)=2n\,t(JX,Y,Z).
\end{aligned}
\]
Next, since $\tau$ is symmetric,
\[
\begin{aligned}
g([\RR(\xi,Z),\tau]X,Y)=g(\RR(\xi,Z)\tau X,Y)
 -g(\RR(\xi,Z)X,\tau Y).
\end{aligned}
\]
Applying \eqref{eq:real-mixedcurvature}
to each term, we find
\[
\begin{aligned}
g(\RR(\xi,Z)\tau X,Y)=t(Z,Y,\tau X)-t(Z,\tau X,Y),
\end{aligned}
\]
and
\[
\begin{aligned}
g(\RR(\xi,Z)X,\tau Y)=t(Z,\tau Y,X)-t(Z,X,\tau Y).
\end{aligned}
\]
Therefore we have
\[
\begin{aligned}
g([\RR(\xi,Z),\tau]X,Y)=t(Z,Y,\tau X)+t(Z,X,\tau Y)-t(Z,\tau X,Y)-t(Z,\tau Y,X).
\end{aligned}
\]
Finally, the definition of $t$ gives
$-g((D_{\tau Z}\tau)X,Y)
=-t(X,Y,\tau Z)$.
Combining these identities in
\eqref{eq:real-commutator}, we obtain
\begin{align}
\label{eq:real-tdot}
(D_\xi t)(X,Y,Z)
={}&2n\,t(JX,Y,Z)+t(Z,Y,\tau X)+t(Z,X,\tau Y)\notag\\
&-t(Z,\tau X,Y)-t(Z,\tau Y,X)-t(X,Y,\tau Z).
\end{align}
Replacing $(X,Z)$ by $(JX,JZ)$ in \eqref{eq:real-tdot} and using $J^2=-I$ and $\tau J=-J\tau$ on $\HH$, we get
\begin{equation*}\tag{\ref{eq:real-tdot}$'$}\label{eq:real-tdotJ}
\begin{aligned}
(D_\xi t)(JX,Y,JZ)
={}&-2n\,t(X,Y,JZ)-t(JZ,Y,J\tau X)
  +t(JZ,JX,\tau Y)\\
&+t(JZ,J\tau X,Y)
  -t(JZ,\tau Y,JX)+t(JX,Y,J\tau Z).
\end{aligned}
\end{equation*}
Since $DJ=0$, we have $(D_\xi\bB)(X,Y,Z)=(D_\xi t)(X,Y,Z)-(D_\xi t)(JX,Y,JZ)$. We subtract \eqref{eq:real-tdotJ} from \eqref{eq:real-tdot}, recalling that $t(JU,JV,W)=-t(U,V,W)$ by \eqref{eq:t-symmetries}. The terms containing $2n$ combine into $2n\,\bB(JX,Y,Z)$. Moreover,
\begin{align*}
 t(Z,Y,\tau X)+t(JZ,Y,J\tau X)&=\bT(Z,Y,X),\\
 -t(X,Y,\tau Z)-t(JX,Y,J\tau Z)&=-\bT(X,Y,Z),\\
 -t(Z,\tau X,Y)-t(JZ,J\tau X,Y)&=0,\\
 t(Z,X,\tau Y)-t(JZ,JX,\tau Y)&=2t(Z,X,\tau Y),\\
 -t(Z,\tau Y,X)+t(JZ,\tau Y,JX)&=-\bB(Z,\tau Y,X)=-\bB(Z,X,\tau Y),
\end{align*}
where the last equality is the total symmetry of $\bB$ (Lemma~\ref{lem:real-Codazzi}). The last two lines add up to $\bC(Z,X,\tau Y)=\bT(Z,X,Y)$, because $2t-\bB=\bC$. Since $\bT$ is symmetric in its first two arguments, we conclude that
\begin{equation}\label{eq:real-Bdot}
 \begin{split}
 (D_\xi\bB)(X,Y,Z)={}&2n\,\bB(JX,Y,Z)+\bT(X,Z,Y)+\bT(Y,Z,X)-\bT(X,Y,Z).
 \end{split}
\end{equation}
Since $\bB$ is totally symmetric, so is
$D_\xi\bB$. Moreover,
$\bB(JX,Y,Z)=\bB(JX,Z,Y)$.
Interchanging $Y$ and $Z$ in
\eqref{eq:real-Bdot} and subtracting,
we therefore obtain
\[
\begin{aligned}
0={}&\bT(X,Z,Y)-\bT(X,Y,Z)+\bT(Y,Z,X)-\bT(Z,Y,X)-\bT(X,Y,Z)+\bT(X,Z,Y).
\end{aligned}
\]
Since $\bT$ is symmetric in its first
two arguments, the above equality becomes
\[
2\bigl(\bT(X,Z,Y)-\bT(X,Y,Z)\bigr)=0.
\]
Thus $\bT$ is also symmetric in its last
two arguments, and consequently it is
totally symmetric.
\end{proof}

\begin{lemma}\label{lem:pointwise}
Let
\begin{equation}\label{eq:real-Fdef}
 \bF(X,Y,Z)=g((D_ZH)X,Y)-g((D_{JZ}H)JX,Y),
\end{equation}
for $X,Y,Z\in\HH$.
Then
\begin{equation}\label{eq:keyH}
 \bF(X,Y,HZ)=\bF(X,HY,Z).
\end{equation}
\end{lemma}

\begin{proof}

Since $DJ=0$, the definition of $\bC$ gives
\[
(D_\xi\bC)(X,Y,Z)
=(D_\xi t)(X,Y,Z)+(D_\xi t)(JX,Y,JZ).
\]
Adding \eqref{eq:real-tdotJ} to \eqref{eq:real-tdot}, the terms containing
$2n$ combine into $2n\,\bC(JX,Y,Z)$,
while the last terms give
$-\bB(X,Y,\tau Z)$.
To simplify the remaining terms, we first
observe that, by \eqref{eq:t-symmetries},
$t(JU,JV,W)=-t(U,V,W)$.
Hence
\begin{gather*}
t(Z,X,\tau Y)+t(JZ,JX,\tau Y)=0,\\
-t(Z,\tau X,Y)+t(JZ,J\tau X,Y)=-2t(Z,\tau X,Y).
\end{gather*}
The other terms combine into
\[
\bB(Z,Y,\tau X)-\bC(Z,\tau Y,X).
\]
Since $\bB$ is totally symmetric by
Lemma~\ref{lem:real-Codazzi}, we have
\[
\begin{aligned}
&\bB(Z,Y,\tau X)-2t(Z,\tau X,Y)=\bB(Z,\tau X,Y)-2t(Z,\tau X,Y)=-\bC(Z,\tau X,Y).
\end{aligned}
\]
Collecting all contributions yields
\begin{equation}\label{eq:real-Cdot}
\begin{split}
(D_\xi\bC)(X,Y,Z)
={}&2n\,\bC(JX,Y,Z)-\bC(Z,\tau Y,X)-\bC(Z,\tau X,Y)-\bB(X,Y,\tau Z).
\end{split}
\end{equation}
Differentiating
the equality \eqref{eq:real-T}, and using
\eqref{eq:C-Jslots}, \eqref{eq:real-Cdot} and
$D_\xi\tau=-2nJ\tau$, we obtain
\begin{equation}\label{eq:real-Tdot}
 \begin{split}
 (D_\xi\bT)(X,Y,Z)={}&4n\,\bT(JX,Y,Z)-\bB(X,Y,HZ)\\
 &-\bC(\tau Z,\tau Y,X)-\bC(\tau Z,\tau X,Y).
 \end{split}
\end{equation}
By Lemma~\ref{lem:secondcodazzi}, $D_\xi\bT$ is totally symmetric.
The same is true for $\bT(JX,Y,Z)$ by \eqref{eq:T-Jslots}.
Interchanging $Y,Z$ in \eqref{eq:real-Tdot} therefore gives
\begin{equation}\label{eq:real-secondconstraint}
 \begin{split}
 0={}&\bB(X,Y,HZ)-\bB(X,Z,HY)+\bC(\tau X,\tau Z,Y)-\bC(\tau X,\tau Y,Z).
 \end{split}
\end{equation}
Replacing $Y$ by $\tau Y$ and using the symmetry of $\bT$, we get
\begin{align*}
 \bC(\tau X,\tau Z,\tau Y)=\bT(\tau X,\tau Z,Y)=\bT(\tau X,Y,\tau Z)=\bC(\tau X,Y,HZ).
\end{align*}
Furthermore, the symmetry of $\bB$ and commutation of $H$ with $\tau$ give
$\bB(X,Z,H\tau Y)=\bB(X,\tau HY,Z)$. Hence
\eqref{eq:real-secondconstraint} becomes
\begin{align}\label{identita-finale1}
 &\bB(X,\tau Y,HZ)+\bC(\tau X,Y,HZ)=\bB(X,\tau HY,Z)+\bC(\tau X,HY,Z).
\end{align}
To conclude, we express $\bF$ in terms of
$\bB$ and $\bC$. Since $H=\tau^2$, we have
$D_ZH=(D_Z\tau)\circ\tau
     +\tau\circ(D_Z\tau)$.
Using the symmetry of $\tau$, this gives
\[
g((D_ZH)X,Y)
=t(\tau X,Y,Z)+t(X,\tau Y,Z).
\]
Similarly, since $\tau J=-J\tau$,
\[
\begin{aligned}
g((D_{JZ}H)JX,Y)
=-t(J\tau X,Y,JZ)+t(JX,\tau Y,JZ).
\end{aligned}
\]
Subtracting these identities and recalling
the definitions of $\bB$ and $\bC$, we obtain
\begin{align}
\bF(X,Y,Z) \notag
={}&t(\tau X,Y,Z)+t(J\tau X,Y,JZ)+t(X,\tau Y,Z)-t(JX,\tau Y,JZ)\\
={}&\bC(\tau X,Y,Z)+\bB(X,\tau Y,Z). \label{identita-finale2}
\end{align}
Substituting \eqref{identita-finale2} into \eqref{identita-finale1} we obtain \eqref{eq:keyH}.
\end{proof}

\begin{proposition}\label{prop:moments}
For every $m\ge1$ and every horizontal vector field $Z$,
\begin{equation}\label{eq:moment-intrinsic}
 \dd s_m(Z)=m\,\dd s_1(H^{m-1}Z).
\end{equation}
Consequently each function $s_m$ is constant on $M$.
\end{proposition}

\begin{proof}
Fix $p\in M$ and choose an orthonormal eigenbasis
$\{e_A\}_{A\in\{1,\ldots,2n\}}$ of $H_p$, with $H_pe_A=q_Ae_A$.
The endomorphism $D_XH$ is symmetric and commutes with $J$.
Thus $(D_XH)J$ is skew-adjoint, and
\[
 g((D_{JZ}H)Je_A,e_A)=0.
\]
Taking $X=Y=e_A$ in \eqref{eq:keyH} gives, at $p$,
\begin{equation}\label{eq:real-diagonal}
 g((D_{HZ}H)e_A,e_A)=q_Ag((D_ZH)e_A,e_A).
\end{equation}
Then, iterating \eqref{eq:real-diagonal},
\begin{align*}
 \dd s_1(H^{m-1}Z)
 =\sum_A q_A^{m-1}g((D_ZH)e_A,e_A)
 =\tr_{\HH}(H^{m-1}D_ZH)=\frac1m\,\dd s_m(Z).
\end{align*}
By Lemma~\ref{lem:einstein}, $s_1=2n-\lambda$ is constant. Hence every
horizontal derivative of $s_m$ vanishes. Moreover,
\[
 \xi(s_m)=m\tr_{\HH}(H^{m-1}D_\xi H)=0.
\]
Thus $\dd s_m=0$, and since $M$ is supposed to be connected, we conclude that $s_m$ is constant on $M$.
\end{proof}

We now have all the tools for proving Theorem~\ref{thm:spectrum-intro}.

\begin{proof}[Proof of Theorem~\ref{thm:spectrum-intro}]
Write the characteristic polynomial of the real endomorphism
$H|_{\HH}$ as
\[
 p_H(t)=t^{2n}-c_1t^{2n-1}+\cdots+(-1)^{2n}c_{2n}.
\]
By Newton's identities we have
\[
 j c_j=\sum_{m=1}^{j}(-1)^{m-1}c_{j-m}s_m,
 \qquad c_0=1.
\]
Proposition~\ref{prop:moments} implies that every coefficient $c_j$ is
constant. Thus the spectrum of $H$, including multiplicities, is
constant. Since $H=h^2$ and $hJ=-Jh$, a positive eigenvalue $q$ of $H$
corresponds to the two eigenvalues $\sqrt q,-\sqrt q$ of $h$ with equal
multiplicities. The kernel is also of constant dimension. Hence the
spectrum of $h$ on $TM$ is constant.

Let $\sigma(h)$ denote the fixed set of
distinct eigenvalues of $h$. For each
$\alpha\in\sigma(h)$, the orthogonal
projection onto the corresponding eigenspace
is given by
\[
P_\alpha=
\prod_{\substack{\beta\in\sigma(h)\\
                 \beta\ne\alpha}}
\frac{h-\beta\Id}{\alpha-\beta}.
\]
Since the eigenvalues are constant, $P_\alpha$
is a smooth endomorphism of constant rank.
Therefore its image
$E_\alpha=\ker(h-\alpha\Id)$
is a smooth subbundle of $TM$.
\end{proof}

\section{Vanishing of the pseudo-Hermitian torsion}\label{sec:vanishing}

By Theorem~\ref{thm:spectrum-intro}, $h$ has constant rank.
Since its nonzero eigenvalues occur in opposite pairs, this rank is
$2r$ for some $0\le r\le n$. We first exclude $r=n$ without assuming
compactness. We then treat the intermediate cases $0<r<n$, where the compactness assumption is used and, by the examples of Section~\ref{sec:examples}, cannot be dropped.

We first establish two identities for the $1$-form
used in the case $r=n$.

\begin{lemma}\label{lem:alpha-identities}
Let $M^{2n+1}$ be an Einstein contact metric
manifold with integrable associated CR structure,
and let $\mathcal U\subseteq M$ be an open set
on which $\tau|_{\HH}$ is invertible.
Define the $1$-form
\begin{equation}\label{eq:alpha}
\alpha(X)=\frac14\tr_{\HH}
\bigl(J\tau^{-1}D_X\tau\bigr).
\end{equation}
Then, for every horizontal vector field $Y$
on $\mathcal U$,
\begin{equation}\label{eq:alpha-horizontal-trace}
\tr_{\HH}(\tau^{-1}D_Y\tau)+4\alpha(JY)
=2g(\delta\tau,\tau^{-1}Y).
\end{equation}
Moreover,
\begin{equation}\label{eq:dalpha-real}
\dd\alpha=-\rho
\end{equation}
on $\mathcal U$.
\end{lemma}
\begin{proof}
Fix $p\in\mathcal U$. Since $\tau$ is symmetric, we can
choose an orthonormal basis $\{e_A\}_{A=1,\ldots,2n}$
of $\HH_p$ such that
\[
\tau e_A=\mu_Ae_A,\qquad \mu_A\ne0.
\]
By the definition of $\alpha$ we have, at $p$,
\[
\begin{aligned}
\tr_{\HH}(\tau^{-1}D_Y\tau)+4\alpha(JY)
&=
\sum_A\bigl[
t(\tau^{-1}e_A,e_A,Y)
+t(J\tau^{-1}e_A,e_A,JY)
\bigr]\\
&=\sum_A\bC(\tau^{-1}e_A,e_A,Y).
\end{aligned}
\]
The definition and total symmetry of $\bT$ (see Lemma~\ref{lem:secondcodazzi}) give
\begin{align}
\sum_A\bC(\tau^{-1}e_A,e_A,Y)
&=\sum_A\bT(\tau^{-1}e_A,e_A,\tau^{-1}Y)
\notag\\
&=\sum_A\bT(\tau^{-1}e_A,\tau^{-1}Y,e_A)
\notag\\
&=\sum_A\bC(\tau^{-1}e_A,\tau^{-1}Y,\tau e_A).
\label{eq:alpha-symmetry}
\end{align}
Expanding $\bC$ and using
$\tau e_A=\mu_Ae_A$, we obtain
\[
\begin{aligned}
\sum_A\bC(\tau^{-1}e_A,\tau^{-1}Y,\tau e_A)&\quad=
\sum_A\bigl[
t(\tau^{-1}e_A,\tau^{-1}Y,\tau e_A)
+t(J\tau^{-1}e_A,\tau^{-1}Y,J\tau e_A)
\bigr]\\
&\quad=
\sum_A\bigl[
t(e_A,\tau^{-1}Y,e_A)
+t(Je_A,\tau^{-1}Y,Je_A)
\bigr]\\
&\quad=
g\left(
\sum_A(D_{e_A}\tau)e_A
+\sum_A(D_{Je_A}\tau)Je_A,\,
\tau^{-1}Y
\right)\\
&\quad=2g(\delta\tau,\tau^{-1}Y),
\end{aligned}
\]
where we used the facts that the factors $\mu_A$ and $\mu_A^{-1}$
cancel by the multilinearity of $t$, and
$\{Je_A\}$ is also an orthonormal basis of
$\HH_p$. Since $p$ is arbitrary, this proves
\eqref{eq:alpha-horizontal-trace}.
We next prove \eqref{eq:dalpha-real}. Put
\[
G_X=\tau^{-1}D_X\tau.
\]
Differentiating $\tau J=-J\tau$ and using
$DJ=0$, we find that $D_X\tau$ also
anticommutes with $J$. Consequently,
\[
G_XJ=JG_X.
\]
Differentiating $\tau^{-1}\tau=I$ gives
\[
D_X(\tau^{-1})
=-\tau^{-1}(D_X\tau)\tau^{-1}.
\]
Using this identity and the definition of
curvature, we obtain
\[
\begin{aligned}
D_X(G_Y)-D_Y(G_X)-G_{[X,Y]}&=
\tau^{-1}\bigl(
D_X(D_Y\tau)-D_Y(D_X\tau)-D_{[X,Y]}\tau
\bigr)-[G_X,G_Y]\\
&=\tau^{-1}[\RR(X,Y),\tau]-[G_X,G_Y].
\end{aligned}
\]
Since $DJ=0$, the definition of the exterior
derivative of $\alpha$ now yields

\begin{equation}\label{alpha-trace}
4\dd\alpha(X,Y)
=\tr_{\HH}
\bigl(J\tau^{-1}[\RR(X,Y),\tau]\bigr)-\tr_{\HH}\bigl(J[G_X,G_Y]\bigr).
\end{equation}
Since $G_Y$ commutes with $J$, we have
$\tr_{\HH}\bigl(J[G_X,G_Y]\bigr)
=\tr_{\HH}\bigl([JG_X,G_Y]\bigr)=0$,
because the trace of a commutator vanishes.
As for the first term on the right-hand side of \eqref{alpha-trace}, using   $\tau J\tau^{-1}=-J$, we have
\[
\begin{aligned}
\tr_{\HH}
\bigl(J\tau^{-1}[\RR(X,Y),\tau]\bigr)&=
\tr_{\HH}\bigl(\tau J\tau^{-1}\RR(X,Y)\bigr)
-\tr_{\HH}\bigl(J\RR(X,Y)\bigr)=-2\tr_{\HH}\bigl(J\RR(X,Y)\bigr).
\end{aligned}
\]
Thus, by \eqref{eq:rho-def},
\[
4\dd\alpha(X,Y)
=-2\tr_{\HH}\bigl(J\RR(X,Y)\bigr)
=-4\rho(X,Y),
\]
which proves \eqref{eq:dalpha-real}.
\end{proof}

\begin{proposition}\label{prop:fullrank}
Let $M^{2n+1}$ be an Einstein contact metric
manifold with integrable associated CR structure,
where $n\ge2$. Then $h|_{\HH}$ is not invertible
at any point.
\end{proposition}

\begin{proof}
Suppose that $h|_{\HH}$ is invertible at some
point. By Theorem~\ref{thm:spectrum-intro},
it is invertible everywhere. The same holds
for $\tau|_{\HH}$, because $\tau^2=h^2$.
We may therefore
apply Lemma~\ref{lem:alpha-identities} on $M$
and consider the globally defined one-form
$\alpha$ in \eqref{eq:alpha}.
We first show that
$\tr_{\HH}(\tau^{-1}D_Y\tau)=0$ for any vector field $Y$ on $M$.
Indeed, since $H=\tau^2$, we have $D_YH=(D_Y\tau)\tau+\tau(D_Y\tau)$, and hence
\[
2\tr_{\HH}(\tau^{-1}D_Y\tau)=\tr_{\HH}(H^{-1}D_YH)=Y\bigl(\log\det(H|_{\HH})\bigr)=0,
\]
because $\det(H|_{\HH})$ is a positive constant by Theorem~\ref{thm:spectrum-intro}.
By \eqref{eq:divtau}, we also have
$\delta\tau=0$. Substituting these
vanishing identities into
\eqref{eq:alpha-horizontal-trace},
we obtain
$4\alpha(JY)=0$
for any $Y\in\Gamma(\HH)$.
Since $J^2=-I$ on $\HH$, it follows that
\begin{equation}\label{vanishing2}
\alpha|_{\HH}=0.
\end{equation}
On the other hand,
$D_\xi\tau=-2nJ\tau$ and
$\tau^{-1}J\tau=-J$ imply
\begin{equation}\label{eq:alpha-xi}
\alpha(\xi)
=\frac14\tr_{\HH}
  \bigl(J\tau^{-1}(-2nJ\tau)\bigr)=-\frac n2\tr_{\HH}I=-n^2.
\end{equation}
By \eqref{vanishing2} and \eqref{eq:alpha-xi}, we conclude that
\begin{equation}\label{eq:alpha-eta}
\alpha=-n^2\eta,\qquad
\dd\alpha=-2n^2\omega.
\end{equation}
Thus,  \eqref{eq:dalpha-real}, \eqref{eq:alpha-eta} and 
\eqref{eq:webster-einstein} imply
\[
(\lambda+2)\omega(X,Y)
=\rho(X,Y)
=-\dd\alpha(X,Y)
=2n^2\omega(X,Y)
\]
for all horizontal vector fields $X,Y$.
Since $\omega$ is nondegenerate on $\HH$,
we conclude that
\begin{equation}\label{lamnda1}
\lambda=2n^2-2.
\end{equation}
However, \eqref{eq:ricxi} and the Einstein
condition give
\begin{equation}\label{lamnda2}
\lambda=2n-\norm{h}^2\le2n,
\end{equation}
Comparing \eqref{lamnda1} with \eqref{lamnda2} we find a contradiction, since
 $2n^2-2>2n$ for $n\ge2$.
\end{proof}

Suppose now that $\rank h=2r$ with $0<r<n$, and put $k=n-r$.
By Theorem~\ref{thm:spectrum-intro}, the real bundles
\begin{equation}\label{eq:real-splitting}
 \EE=\im\tau=\im h,\qquad
 \KK=\ker\tau\cap\HH=\ker h\cap\HH
\end{equation}
are smooth, orthogonal and $J$-invariant, with ranks $2r$ and $2k$.
Set $S=\tau|_{\EE}$. $S$  is a symmetric invertible endomorphism
anticommuting with $J|_{\EE}$, and its eigenvalues are constant.

Let $D^{\EE}$ and $D^{\KK}$ be the connections induced by orthogonal
projection of $D$, and put $D^0=D^{\EE}\oplus D^{\KK}$ on $\HH$.
We use the same notation for the induced connections on tensor
bundles constructed from $\EE$ and $\KK$. Define the second fundamental
form of $\EE$ by
\begin{equation}\label{eq:real-secondfundamental}
 \BB_XY=\pi_{\KK}D_XY,
\end{equation}
for $X\in\Gamma( TM)$ and $Y\in\Gamma(\EE)$.
Since $D$ is a metric connection and
$\EE$ and $\KK$ are orthogonal, we have
\[
\pi_{\EE}D_XZ=-\BB_X^*Z,
\]
for all $Z\in\Gamma(\KK)$,
where $\BB_X^*:\KK\longrightarrow\EE$ denotes the
metric adjoint of $\BB_X$.
Moreover, since $DJ=0$, one also has $\BB_X(JY)=J\BB_XY$.

\begin{lemma}\label{lem:Bcomponents}
For all $X,Y\in\Gamma(\EE)$ and $Z\in\Gamma(\KK)$,
\begin{equation}\label{eq:Bzeros}
 \BB_\xi=0,\qquad \BB_{JX}Y=-J\BB_XY,\qquad
 \BB_{JZ}Y=J\BB_ZY.
\end{equation}
Consequently, for any $Y\in\Gamma(\EE)$ the  $(1,1)$-tensor field  defined by
\begin{equation}\label{eq:real-Ldef}
 \LL_YZ=\BB_ZY,
\end{equation}
for any $Z\in\Gamma(\KK)$,
satisfies
\begin{equation}\label{eq:L-J}
 \LL_YJ=J\LL_Y,\qquad \LL_{JY}=J\LL_Y.
\end{equation}
\end{lemma}

\begin{proof}
The identity $D_\xi H=0$ implies that $D_\xi$ preserves $\KK$ and
its orthogonal complement, proving $\BB_\xi=0$.
If $Z,W\in\Gamma(\KK)$, then $\tau Z=\tau W=0$.
Hence, for any horizontal vector field $X$,
\[
\begin{aligned}
t(Z,W,X)
=g((D_X\tau)Z,W)=g(D_X(\tau Z)-\tau(D_XZ),W)=-g(D_XZ,\tau W)=0,
\end{aligned}
\]
where the last equality follows from the
symmetry of $\tau$.
 Since $JZ\in\Gamma(\KK)$ as well, this yields $\bB(Z,W,Y)=0$ for every horizontal $Y$. The symmetry of $\bB$ therefore
gives $\bB(Y,Z,W)=0$ for $Y\in\EE$. Its expansion is
\[
 g\bigl(\BB_W(\tau Y)+J\BB_{JW}(\tau Y),Z\bigr)=0.
\]
Since $Z$ is arbitrary and $S$ is invertible, this implies
$\BB_{JW}Y=J\BB_WY$.
Next, $\bT(X,Y,Z)=0$ when $Z\in\KK$. By the symmetry of $\bT$,
$\bC(X,Z,\tau Y)=0$ for $X,Y\in\EE$ and $Z\in\KK$.
Expanding this equality gives
\[
 g\bigl(\BB_{\tau Y}(\tau X)-J\BB_{J\tau Y}(\tau X),Z\bigr)=0.
\]
Invertibility of $S$ now yields $\BB_{JX}Y=-J\BB_XY$.
The first equality in \eqref{eq:L-J} follows from the last equality
in \eqref{eq:Bzeros}, and the second follows from
$\BB_Z(JY)=J\BB_ZY$.
\end{proof}

We regard $\LL$ as a bundle map
$\EE\longrightarrow\operatorname{End}(\KK)$. Its covariant derivative
is defined by
\begin{equation}\label{eq:L-derivative}
 \begin{split}
 (D_X^0\LL)(Y)Z={}&D_X^{\KK}(\LL_YZ)
            -\LL_Y(D_X^{\KK}Z)-\LL_{D_X^{\EE}Y}Z,
 \end{split}
\end{equation}
for $Y\in\Gamma(\EE)$ and $Z\in\Gamma(\KK)$.

\begin{lemma}\label{lem:Cequations}
For any $X,Y\in\Gamma(\EE)$, the endomorphisms $\LL_Y$ satisfy
\begin{align}
 (D_\xi^0\LL)(Y)&=-\LL_{SY}^{*},\label{eq:Lxi}\\
 (D_X^0\LL)(Y)-J(D_{JX}^0\LL)(Y)+2\LL_Y\LL_X
  &=-2g(SY,JX)\Id_{\KK}-2g(SY,X)J|_{\KK},
  \label{eq:Lhorizontal}
\end{align}
where $\LL_Y^*$ denotes the metric adjoint of the endomorphism
$\LL_Y:\KK\longrightarrow\KK$.
\end{lemma}

\begin{proof}
Fix $p\in M$ and choose local orthonormal
frames $\{e_A\}$ of $\EE$ and $\{v_\mu\}$ of
$\KK$ such that
\[
(D^0_Xe_A)_p=0,\qquad
(D^0_Xv_\mu)_p=0
\]
for every $X\in T_pM$.
Since the identities to be proved are
tensorial, we may extend arbitrary vectors
$X_p,Y_p\in\EE_p$ and $Z_p,W_p\in\KK_p$
to local sections $X,Y,Z,W$ having constant
coefficients in these frames. Thus
\[
(D^0_VX)_p=(D^0_VY)_p
=(D^0_VZ)_p=(D^0_VW)_p=0
\]
for every $V\in T_pM$.
For $Y\in\Gamma(\EE)$ and $Z,W\in\Gamma(\KK)$, \eqref{eq:real-mixedcurvature} gives
\begin{equation}\label{eq:Lxi-curvature}
 g(\RR(\xi,Z)Y,W)=-g(\LL_{SY}W,Z).
\end{equation}
Indeed, $t(Z,W,Y)=0$, while
$g((D_W\tau)Z,Y)=g(\BB_W(SY),Z)$.
Since $\BB_\xi=0$, $D_\xi$ preserves the splitting.
Moreover $T^D(\xi,Z)=0$ and $D\xi=0$, so $[\xi,Z]_p=0$.
The curvature definition, projected to $\KK$, therefore gives
\[
 \pi_{\KK}\RR(\xi,Z)Y=(D_\xi^0\LL)(Y)Z
\]
at $p$. Comparing this with \eqref{eq:Lxi-curvature} proves
\eqref{eq:Lxi}.
For the second identity, $T^D(X,Z)=0$ when $X\in\Gamma(\EE)$, $Z\in\Gamma(\KK)$.
At the point $p$, the normality of the frames gives
\[
 D_XZ=-\BB_X^*Z,\qquad D_ZX=\LL_XZ,
 \qquad [X,Z]=-\BB_X^*Z-\LL_XZ.
\]
The first derivative here may be nonzero: normality for $D^0$ does not
annihilate the second fundamental form. Projecting the definition of
$\RR(X,Z)Y$ to $\KK$ gives
\begin{equation}\label{eq:L-curvature-projection}
 \pi_{\KK}\RR(X,Z)Y=(D_X^0\LL)(Y)Z-(D_Z^0\BB)(X,Y)+\BB_{\BB_X^*Z}Y+\LL_Y\LL_XZ.
\end{equation}
Here $D_Z^0\BB$ differentiates the restriction of $\BB$ to
$\EE\times\EE$ and uses $D^{\KK}$ on its values.
Consider the following combination of curvature terms:
\[
 \mathcal Q(X,Y,Z)=\pi_{\KK}\bigl(
 \RR(X,Z)Y-\RR(JX,JZ)Y
          -J\RR(JX,Z)Y-J\RR(X,JZ)Y\bigr).
\]
Substitute \eqref{eq:L-curvature-projection} into each term.
The derivatives of $\BB|_{\EE\times\EE}$ cancel by
$\BB_{JX}Y=-J\BB_XY$ and $D^0J=0$.
The terms quadratic in this part of $\BB$ cancel as well, using
$\BB_{JX}^{*}=\BB_X^{*}J$.
The other terms give
\begin{equation}\label{eq:L-comparison-derivatives}
 \mathcal Q(X,Y,Z)
 =2\bigl((D_X^0\LL)(Y)-J(D_{JX}^0\LL)(Y)
                                      +2\LL_Y\LL_X\bigr)Z.
\end{equation}
Alternatively, decompose $\RR=\RR_0+\RR_\tau$ as in
Lemma~\ref{lem:curvature-decomposition}. The contribution of $\RR_0$
is zero by \eqref{eq:R0-J}. Since $\tau Z=0$, direct substitution of
\eqref{eq:curvature-model} gives
\begin{equation}\label{eq:L-comparison-curvature}
 \mathcal Q(X,Y,Z)=-4g(SY,JX)Z-4g(SY,X)JZ.
\end{equation}
Comparison of \eqref{eq:L-comparison-derivatives} and
\eqref{eq:L-comparison-curvature}, for all $Z\in\KK_p$, proves
\eqref{eq:Lhorizontal}. Since $p$ is arbitrary, both identities hold
on $M$.
\end{proof}

Define the  $1$-form $\nu$ on $\EE$ and the  vector field $W$ by
\begin{equation}\label{eq:real-nuW}
 \nu(Y)=\tr_{\KK}\LL_Y,\qquad W=\frac12S^{-1}\nu^\sharp.
\end{equation}
The metric dual $\nu^\sharp$ is taken in $\EE$. For a real orthonormal
frame $\{e_A\}_{A\in\{1,\ldots,2r\}}$ of $\EE$, put
\begin{align}
 U&=\frac12\sum_A\tr_{\KK}(\LL_{e_A}^*\LL_{e_A})\ge0,
          \label{eq:globalU}\\
 f_0&=\frac12\sum_A\tr_{\KK}(\LL_{S^{-1}e_A}\LL_{e_A}),
          \label{eq:globalf0}\\
 f_1&=-\frac12\sum_A\tr_{\KK}(J\LL_{S^{-1}e_A}\LL_{e_A}).
          \label{eq:globalf1}
\end{align}
These expressions use the metrics of $\EE$ and $\KK$ and are
independent of the orthonormal frame. Thus $U,f_0,f_1$ are
global  functions and $W$ is a global vector field.
The endomorphism $S$ has constant nonzero eigenvalues, hence
$|\det S|$ is a positive constant.

\Needspace{14\baselineskip}
\begin{lemma}\label{lem:inverse-divergence}
One has
\begin{equation}\label{eq:real-inverse-div}
 \sum_A(D_{e_A}^0S^{-1})e_A=-S^{-1}\nu^\sharp.
\end{equation}
Hence,
\begin{equation}\label{eq:divJW}
 \divg(-JW)=-\frac12\sum_A(D_{e_A}^0\nu)(JS^{-1}e_A).
\end{equation}
\end{lemma}

\begin{proof}
For an orthonormal frame $\{v_\mu\}_{\mu\in\{1,\ldots,2k\}}$ of $\KK$, the
$\EE$ component of $\delta\tau=0$ reads
\begin{equation}\label{eq:deltaES}
 \delta_{\EE}S:=\sum_A(D_{e_A}^0S)e_A=-S\nu^\sharp.
\end{equation}
Indeed,
\[
 \sum_\mu(D_{v_\mu}\tau)v_\mu
       =-S\sum_\mu\pi_{\EE}D_{v_\mu}v_\mu=S\nu^\sharp,
\]
because the scalar product of the middle sum with $Y\in\EE$ equals
$-\sum_\mu g(\BB_{v_\mu}Y,v_\mu)=-\nu(Y)$.
Note that $t(X,Y,Z)=g((D_Z^0S)X,Y)$ for all $X,Y\in\Gamma(\EE)$ and every horizontal vector field $Z$, because $\tau$ is symmetric and preserves $\EE$.
Since $|\det S|$ is constant,
\begin{equation}\label{eq:real-determinant}
 \tr_{\EE}(S^{-1}D_Z^0S)=Z(\log|\det S|)=0.
\end{equation}
Fix $Y\in\Gamma(\EE)$. Since $\bB=2t-\bC$, using \eqref{eq:real-determinant}, the definition and the total symmetry of $\bT$ (Lemma~\ref{lem:secondcodazzi}), and finally \eqref{eq:deltaES}, we get
\begin{align*}
\sum_A\bB(S^{-1}e_A,e_A,Y)
&=2\tr_{\EE}(S^{-1}D_Y^0S)-\sum_A\bC(S^{-1}e_A,e_A,Y)
 =-\sum_A\bT(S^{-1}e_A,e_A,S^{-1}Y)\\
&=-\sum_A\bT(S^{-1}e_A,S^{-1}Y,e_A)
 =-\sum_A\bC(S^{-1}e_A,S^{-1}Y,Se_A)\\
&=-2\sum_A t(e_A,S^{-1}Y,e_A)
 =-2g(\delta_{\EE}S,S^{-1}Y)=2\nu(Y).
\end{align*}
In the fifth equality we used that $\sum_AS^{-1}e_A\otimes Se_A=\sum_Ae_A\otimes e_A$, because $S$ is symmetric, and that the two summands in the definition \eqref{eq:real-C} of $\bC$ give the same contribution, as one sees by replacing the orthonormal frame $\{e_A\}$ by $\{Je_A\}$ and using $JS=-SJ$ and $JS^{-1}=-S^{-1}J$. By the same change of frame and the total symmetry of $\bB$ (Lemma~\ref{lem:real-Codazzi}),
\[
\sum_A\bB(S^{-1}e_A,e_A,Y)=\sum_A\bB(S^{-1}e_A,Y,e_A)=2\sum_At(S^{-1}e_A,Y,e_A).
\]
Comparing the last two formulas, we obtain $\sum_At(S^{-1}e_A,Y,e_A)=\nu(Y)$ for every $Y\in\Gamma(\EE)$, that is,
$\sum_A(D_{e_A}^0S)S^{-1}e_A=\nu^\sharp$.
Differentiating $S^{-1}S=\Id_{\EE}$ gives
\[
 \sum_A(D_{e_A}^0S^{-1})e_A
 =-S^{-1}\sum_A(D_{e_A}^0S)S^{-1}e_A,
\]
which proves \eqref{eq:real-inverse-div}.
For a vector field $Y\in\Gamma(\EE)$, formula \eqref{eq:connection}
shows that $V\mapsto(\nabla_V-D_V)Y$ has zero trace. The traces over
$\EE$, $\KK$ and $\R\xi$ therefore give
\begin{equation}\label{eq:div-adapted}
 \divg Y=\sum_Ag(D^{\EE}_{e_A}Y,e_A)+\nu(Y).
\end{equation}
The $\KK$ contribution is $\tr_{\KK}\LL_Y=\nu(Y)$, and the
vertical contribution is zero. Apply \eqref{eq:div-adapted} to $-JW=-\frac12JS^{-1}\nu^\sharp$. Since $D^0J=0$ and the endomorphisms $JS^{-1}$ and $J(D^0_{e_A}S^{-1})$ of $\EE$ are symmetric, we get
\[
\sum_Ag\bigl(D^{\EE}_{e_A}(-JW),e_A\bigr)
=-\frac12\,g\Bigl(\nu^\sharp,J\sum_A(D^0_{e_A}S^{-1})e_A\Bigr)
 -\frac12\sum_A(D^0_{e_A}\nu)(JS^{-1}e_A).
\]
By \eqref{eq:real-inverse-div}, the first term on the right-hand side equals $\frac12g(\nu^\sharp,JS^{-1}\nu^\sharp)=-\nu(-JW)$, and \eqref{eq:divJW} follows.
\end{proof}

\Needspace{11\baselineskip}
\begin{lemma}\label{lem:two-global}
The functions and vector field defined in
\eqref{eq:real-nuW}--\eqref{eq:globalf1} satisfy
\begin{align}
 \xi(f_0)&=2nf_1-2U,\label{eq:real-Reeb-contraction}\\
 \divg(-JW)&=-2rk-f_1.\label{eq:real-divJW-final}
\end{align}
\end{lemma}
\begin{proof}

Recall that
$\nu(Y)=\tr_{\KK}(\LL_Y)$.
Since $D^0$ is a metric connection preserving
$\EE$ and $\KK$, differentiation gives
\[
(D_X^0\nu)(Y)
=\tr_{\KK}\bigl((D_X^0\LL)(Y)\bigr).
\]
Moreover, differentiating $\LL_{JY}=J\LL_Y$
and using $D^0J=0$, we obtain
\[
(D_X^0\LL)(JY)=J(D_X^0\LL)(Y).
\]
Therefore,
\[
(D_X^0\nu)(JY)
=\tr_{\KK}\bigl(J(D_X^0\LL)(Y)\bigr).
\]
We now multiply \eqref{eq:Lhorizontal}
by $J$ on the left and take the trace
over $\KK$. Since $J^2=-I$ on $\KK$,
the left-hand side becomes
\[
\begin{aligned}
\tr_{\KK}\bigl(J(D_X^0\LL)(Y)\bigr)
-\tr_{\KK}\bigl(J^2(D_{JX}^0\LL)(Y)\bigr)+2\tr_{\KK}(J\LL_Y\LL_X)
&=(D_X^0\nu)(JY)+(D_{JX}^0\nu)(Y)\\
  &\quad+2\tr_{\KK}(J\LL_Y\LL_X).
\end{aligned}
\]
Since $\tr_{\KK}J=0$ and
$\tr_{\KK}J^2=-2k$, the trace of the
right-hand side is
\[
\begin{aligned}
-2g(SY,JX)\tr_{\KK}J
-2g(SY,X)\tr_{\KK}J^2=4k\,g(SY,X).
\end{aligned}
\]
We therefore obtain
\begin{equation}\label{eq:Ltrace1}
\begin{aligned}
(D_X^0\nu)(JY)+(D_{JX}^0\nu)(Y)+2\tr_{\KK}(J\LL_Y\LL_X)
=4k\,g(SY,X).
\end{aligned}
\end{equation}
Let $\{e_A\}_{A\in\{1,\ldots,2r\}}$ be a local
orthonormal frame of $\EE$, and substitute $X=e_A$ and $Y=S^{-1}e_A$ into \eqref{eq:Ltrace1}. Since
$g(e_A,e_A)=1$ and $\dim\EE=2r$,
summing over $A$ gives
\[
\begin{aligned}
\sum_A(D_{e_A}^0\nu)(JS^{-1}e_A)+\sum_A(D_{Je_A}^0\nu)(S^{-1}e_A)+2\sum_A
\tr_{\KK}(J\LL_{S^{-1}e_A}\LL_{e_A})
=8rk.
\end{aligned}
\]
Replacing $\{e_A\}$ by $\{Je_A\}$
in the second sum and using $S^{-1}J=-JS^{-1}$, we obtain
\[
\begin{aligned}
\sum_A(D_{Je_A}^0\nu)(S^{-1}e_A)=\sum_A(D_{e_A}^0\nu)(JS^{-1}e_A).
\end{aligned}
\]
Furthermore, the definition of $f_1$
gives
\[
\sum_A\tr_{\KK}
(J\LL_{S^{-1}e_A}\LL_{e_A})=-2f_1.
\]
Thus
\[
2\sum_A(D_{e_A}^0\nu)(JS^{-1}e_A)
-4f_1=8rk,
\]
which yields
\begin{equation}\label{eq:nu-contraction1}
\sum_A(D_{e_A}^0\nu)(JS^{-1}e_A)
=4rk+2f_1.
\end{equation}
Finally, substituting this identity into \eqref{eq:divJW},
we obtain
\[
\divg(-JW)
=-\frac12\sum_A
(D_{e_A}^0\nu)(JS^{-1}e_A)=-2rk-f_1,
\]
which is precisely
\eqref{eq:real-divJW-final}.
It remains to prove \eqref{eq:real-Reeb-contraction}.
We first recall two identities needed to
differentiate $f_0$. Since $\BB_\xi=0$,
the connection $D_\xi$ preserves both
$\EE$ and $\KK$. In particular, its
restriction to $\EE$ agrees with $D_\xi^0$.
Thus, by Lemma~\ref{lem:einstein},
\begin{equation}\label{identita1}
D_\xi^0S=(D_\xi\tau)|_{\EE}=-2nJS.
\end{equation}
Differentiating $SS^{-1}=\Id_{\EE}$ and
using $S^{-1}J=-JS^{-1}$, we also obtain
\begin{equation}\label{identita2}
D_\xi^0S^{-1}=-S^{-1}(D_\xi^0S)S^{-1}=2nS^{-1}J=-2nJS^{-1}.
\end{equation}
Fix an arbitrary point $p\in M$ and choose
a local orthonormal frame $\{e_A\}$ of $\EE$
such that
$(D_\xi^0e_A)_p=0$.
Recall that
\[
f_0=\frac12\sum_A
\tr_{\KK}(\LL_{S^{-1}e_A}\LL_{e_A}).
\]
Differentiating this expression along $\xi$
at $p$,
using \eqref{identita1} and \eqref{identita2}, together with
\eqref{eq:Lxi} and $\LL_{JY}=J\LL_Y$,
we obtain
\begin{align}\label{identita5}
\xi(f_0)
={}&\frac12\sum_A\tr_{\KK}\Bigl(
\LL_{(D_\xi^0S^{-1})e_A}\LL_{e_A} +(D_\xi^0\LL)(S^{-1}e_A)\LL_{e_A}+\LL_{S^{-1}e_A}(D_\xi^0\LL)(e_A)
\Bigr)\\
={}&-n\sum_A\tr_{\KK}
(J\LL_{S^{-1}e_A}\LL_{e_A})-\frac12\sum_A\tr_{\KK}
(\LL_{e_A}^*\LL_{e_A})-\frac12\sum_A\tr_{\KK}
(\LL_{S^{-1}e_A}\LL_{Se_A}^*). \nonumber
\end{align}
By the definitions of $f_1$ and $U$,
the first two sums equal $2nf_1$
and $-U$, respectively. To simplify
the last sum, observe that the symmetry
of $S$ gives
\begin{equation}\label{identita3}
\sum_A
g(S^{-1}e_A,e_B)\,g(Se_A,e_C)=g(S^{-1}e_B,Se_C)
=\delta_{BC}.
\end{equation}
Hence,
\begin{equation}\label{identita4}
\sum_A\tr_{\KK}
(\LL_{S^{-1}e_A}\LL_{Se_A}^*)
=\sum_B\tr_{\KK}
(\LL_{e_B}\LL_{e_B}^*)
=\sum_B\tr_{\KK}
(\LL_{e_B}^*\LL_{e_B})=2U.
\end{equation}
Substituting \eqref{identita3} and \eqref{identita4} into
\eqref{identita5}, we conclude
that
\[
\xi(f_0)=2nf_1-U-U
=2nf_1-2U.
\]
\end{proof}

\begin{proposition}\label{prop:intermediate}
Let $M^{2n+1}$ be a compact Einstein contact metric manifold with
integrable associated CR structure. Then the rank of $h$ is  either $0$ or $2n$.
\end{proposition}

\begin{proof}

Suppose that $\rank h=2r$ with $0<r<n$,
and put $k=n-r$.
Since $\nabla\xi=J+\tau$ is traceless, we have $\divg\xi=0$.
 Consider the global 
vector field $W$ and the smooth functions
$U$, $f_0$, and $f_1$, introduced in
\eqref{eq:real-nuW}--\eqref{eq:globalf1}, and define
 the  vector field
\begin{equation}\label{eq:negativefield}
 \mathcal V=-JW+\frac{f_0}{2n}\xi.
\end{equation}
By Lemma~\ref{lem:two-global} we have
\begin{equation}\label{eq:negative-div}
 \divg\mathcal V=(-2rk-f_1)+\frac{2nf_1-2U}{2n}
                    =-2r(n-r)-\frac{U}{n}<0.
\end{equation}
Then integrating over $M$ contradicts
the divergence theorem.
\end{proof}

Finally we have the tools for proving Theorem~\ref{thm:main}.

\begin{proof}[Proof of Theorem~\ref{thm:main}]
Theorem~\ref{thm:spectrum-intro} gives a constant even rank $2r$ for
$h$, with $0\le r\le n$. Proposition~\ref{prop:fullrank} excludes
$r=n$, and Proposition~\ref{prop:intermediate} excludes $0<r<n$.
Thus $h=0$, and \eqref{eq:CRlevi} reduces to
\[
 (\nabla_X\phi)Y=g(X,Y)\xi-\eta(Y)X,
\]
which characterizes the given contact metric structure as Sasakian.
Finally, \eqref{eq:ricxi} gives $\lambda=2n$.
\end{proof}

\section{Noncompact examples}\label{sec:examples}

We prove Theorem~\ref{thm:examples-intro} by an explicit construction.
 In what follows, we put $k=n-1\ge1$.
On $\R^{2n+1}$ take coordinates
$(p,q,x_1,y_1,\ldots,x_k,y_k,t)$ and define
\begin{gather}
 w=n(p^2+q^2)-2pq,\qquad F=1+\frac{w}{2k},\qquad
 f=F^{1/(n+1)},\label{eq:ex-functions}\\
 \Delta=2(n^2-1)F,\qquad
 b=\frac{1}{2k(1+f+\cdots+f^n)}.\label{eq:ex-delta}
\end{gather}
The quadratic form $w$ is positive definite,
so $F$, $f$, and $\Delta$ are positive
everywhere. Moreover, since
\[
w=2k(f^{n+1}-1)
  =2k(f-1)(1+f+\cdots+f^n),
\]
we have
$b=\frac{f-1}{w}$
whenever $(p,q)\ne(0,0)$.
The expression defining $b$ in
\eqref{eq:ex-delta} provides a smooth
extension of this quotient to the origin,
where
$b(0,0)=\frac{1}{2k(n+1)}$.
Thus $b$ is smooth on all of $\mathbb R^2$.
 Set
\begin{gather}
 \eta=f^{-1}\left(\dd t+\sum_{j=1}^k(x_j\dd y_j-y_j\dd x_j)
                  -b(p\dd q-q\dd p)\right),\label{eq:ex-eta}\\
 \alpha=\frac q{\sqrt2}\eta+
 \frac{-((n+1)q^2+2n)\dd p+((n+1)pq+2)\dd q}
      {\sqrt2\,\Delta},\label{eq:ex-alpha}\\
 \beta=\frac p{\sqrt2}\eta+
 \frac{-((n+1)pq+2)\dd p+((n+1)p^2+2n)\dd q}
      {\sqrt2\,\Delta},\label{eq:ex-beta}\\
 \chi_j=f^{-1/2}\dd x_j,\qquad
 \psi_j=f^{-1/2}\dd y_j.\label{eq:ex-chi}
\end{gather}
We now verify that the $2n+1$ differential
forms defined above constitute a global
coframe on $\mathbb R^{2n+1}$.
To this aim, we compute their exterior
product. Since
 $\eta\wedge\eta=0$, we have
\[
\begin{aligned}
\eta\wedge\alpha\wedge\beta
&=\eta\wedge
\left(\alpha-\frac q{\sqrt2}\eta\right)
\wedge
\left(\beta-\frac p{\sqrt2}\eta\right)\\
&=\frac{
((n+1)pq+2)^2
-((n+1)q^2+2n)((n+1)p^2+2n)}
{2\Delta^2}\,
\eta\wedge dp\wedge dq\\
&=-\frac1\Delta\,
\eta\wedge dp\wedge dq,
\end{aligned}
\]
where the last equality follows from 
$\Delta=2(n^2-1)+(n+1)w$.
On the other hand, the definitions of
$\chi_j$ and $\psi_j$ give
\[
\bigwedge_{j=1}^k
(\chi_j\wedge\psi_j)
=f^{-k}\bigwedge_{j=1}^k
(dx_j\wedge dy_j).
\]
Consequently,
\[
\begin{aligned}
\eta\wedge\alpha\wedge\beta
 \wedge\bigwedge_{j=1}^k
 (\chi_j\wedge\psi_j)&=-\frac{f^{-k}}{\Delta}\,
\eta\wedge dp\wedge dq
\wedge\bigwedge_{j=1}^k
(dx_j\wedge dy_j)\\
&\quad=-\frac{f^{-(k+1)}}{\Delta}\,
dt\wedge dp\wedge dq
\wedge\bigwedge_{j=1}^k
(dx_j\wedge dy_j),
\end{aligned}
\]
where, in the last step, we used the definition
of $\eta$ and the fact that all its terms other than
$f^{-1}dt$ vanish when wedged with
$dp\wedge dq$ and the forms
$dx_j\wedge dy_j$.
Since $f>0$ and $\Delta>0$ everywhere,
this exterior product never vanishes.
Thus
$\left\{\eta,\alpha,\beta,\chi_1,\psi_1,
\ldots,\chi_k,\psi_k\right\}$
is a global coframe on $\mathbb R^{2n+1}$.
We then define a Riemannian metric $g$ by
declaring this coframe to be orthonormal, that is
\begin{equation}\label{eq:ex-metric}
g=\eta^2+\alpha^2+\beta^2+
\sum_{j=1}^k(\chi_j^2+\psi_j^2).
\end{equation}
Let $\left\{\xi,E_1,E_2,X_1,Y_1,\ldots,X_k,Y_k\right\}$ be the dual frame and define
\begin{equation}\label{eq:ex-phi}
 \phi\xi=0,\quad \phi E_1=-E_2,\quad \phi E_2=E_1,
 \quad \phi X_j=-Y_j,\quad \phi Y_j=X_j.
\end{equation}
As in Section \ref{sec:prelim}, we set $J=-\phi$. Differentiation gives
\begin{align}
 \dd\eta&=2\alpha\wedge\beta+2\Omega,\label{eq:ex-structure0}\\
 \dd\alpha&=\eta\wedge\alpha-n\eta\wedge\beta
       -\frac{kq}{\sqrt2}\alpha\wedge\beta+\sqrt2q\,\Omega,
       \label{eq:ex-structure1}\\
 \dd\beta&=n\eta\wedge\alpha-\eta\wedge\beta
       -\frac{kp}{\sqrt2}\alpha\wedge\beta+\sqrt2p\,\Omega,
       \label{eq:ex-structure2}\\
 \dd\chi_j&=\frac{p\alpha-q\beta}{\sqrt2}\wedge\chi_j,
 \qquad
 \dd\psi_j=\frac{p\alpha-q\beta}{\sqrt2}\wedge\psi_j,
       \label{eq:ex-structureK}
\end{align}
where we put  $\Omega=\sum_j\chi_j\wedge\psi_j$.
To compute the derivatives of functions
depending on $p$ and $q$, we shall use
the following expressions for $\dd p$
and $\dd q$ in terms of the orthonormal
coframe introduced above:
\begin{align}
 \dd p&=(-p+nq)\eta-\frac{(n+1)p^2+2n}{\sqrt2}\alpha
                       +\frac{(n+1)pq+2}{\sqrt2}\beta,
                       \label{eq:ex-dp}\\
 \dd q&=(-np+q)\eta-\frac{(n+1)pq+2}{\sqrt2}\alpha
                       +\frac{(n+1)q^2+2n}{\sqrt2}\beta.
                       \label{eq:ex-dq}
\end{align}

Equation \eqref{eq:ex-structure0} gives
\[
 \eta\wedge(\dd\eta)^n
 =2^nn!\,\eta\wedge\alpha\wedge\beta
          \wedge\chi_1\wedge\psi_1\wedge\cdots
          \wedge\chi_k\wedge\psi_k\ne0.
\]
Thus $\eta$ is a contact form and 
the Reeb vector field is the vector dual to $\eta$.  More explicitly,
\begin{equation}\label{eq:ex-reeb}
 \xi=\partial_t+(-p+nq)\partial_p+(-np+q)\partial_q.
\end{equation}
Straightforward computations show that \eqref{eq:contactalg}--\eqref{eq:contactnormalization} are satisfied, so that $(\phi,\xi,\eta,g)$ is a contact metric structure on $\mathbb{R}^{2n+1}$. Note that 
the brackets involving $\xi$ are given by
\[
 [\xi,E_1]=-E_1-nE_2,\qquad [\xi,E_2]=nE_1+E_2,
 \qquad [\xi,X_j]=[\xi,Y_j]=0.
\]
Consequently
\begin{equation}\label{eq:ex-h}
 hE_1=-E_2,\qquad hE_2=-E_1,
 \qquad hX_j=hY_j=h\xi=0.
\end{equation}
Thus  $h$ has rank $2$, and eigenvalues $1$ and $-1$,
each of multiplicity one, and $0$ of
multiplicity $2n-1$. In particular, the 
structure is not Sasakian.

We now verify the integrability of the
associated CR structure. Using the structure
equations
\eqref{eq:ex-structure0}--\eqref{eq:ex-structureK}
and Koszul's formula, we compute the
covariant derivative of $\phi$.
Since $\nabla\phi$ is tensorial, it suffices
to perform the calculation on the elements
of the orthonormal frame. For example,
Koszul's formula gives
\begin{equation*}
\nabla_{X_j}X_j
=\frac{-pE_1+qE_2}{\sqrt2},\quad
\nabla_{X_j}Y_j
=-\xi-\frac{qE_1+pE_2}{\sqrt2}.
\end{equation*}
Since $\phi X_j=-Y_j$ and $\phi=-J$
on $\HH$, we obtain
\[
\begin{aligned}
(\nabla_{X_j}\phi)X_j=-\nabla_{X_j}Y_j
  +J\nabla_{X_j}X_j=\xi.
\end{aligned}
\]
The corresponding calculations for the
remaining frame vectors yield, for
$A,B\in\operatorname{span}\{E_1,E_2\}$
and $V,W\in\KK_0$, where $\KK_0=\operatorname{span}\{X_1,Y_1,\ldots,X_k,Y_k\}=\ker h\cap\HH$,
\begin{equation}\label{eq:ex-real-integrability}
(\nabla_A\phi)B=g(A+hA,B)\xi, \quad
(\nabla_A\phi)V=(\nabla_V\phi)A=0,\quad
(\nabla_V\phi)W=g(V,W)\xi.
\end{equation}
By the tensoriality of $\nabla\phi$, these
identities hold for arbitrary smooth sections
$A,B$ of $\operatorname{span}\{E_1,E_2\}$
and $V,W$ of $\KK_0$.
Since $h$ vanishes on $\KK_0$, these
identities imply
\begin{equation}\label{identita7}
(\nabla_U\phi)V=g(U+hU,V)\xi
\end{equation}
for all $U,V\in\Gamma(\HH)$.
The same calculation gives
\begin{equation}\label{identita8}
\nabla_\xi\phi=0,\qquad
(\nabla_U\phi)\xi=-(U+hU),
\end{equation}
for all $U\in\Gamma(\HH)$.
Together, the identities \eqref{identita7} and \eqref{identita8} establish
\eqref{eq:CRlevi}.

Finally, we  verify that the metric $g$ is Einstein.
Let $\{e_A\}$ denote the global orthonormal
frame dual to the coframe constructed above,
and write
\[
[e_A,e_B]=\sum_C c_{AB}^{C}e_C.
\]
The coefficients $c_{AB}^{C}$ are determined
by the structure equations
\eqref{eq:ex-structure0}--\eqref{eq:ex-structureK}.
Koszul's formula gives
\[
2g(\nabla_{e_A}e_B,e_C)
=c_{AB}^{C}-c_{BC}^{A}+c_{CA}^{B}.
\]
We use these identities to compute the
Riemannian curvature tensor, with the
derivatives of $p$ and $q$ given by
\eqref{eq:ex-dp}--\eqref{eq:ex-dq}.

By a direct computation we have that the sectional curvatures
needed to compute the Ricci tensor are those
listed in the following table, where  $s=p^2+q^2$.
\medskip
\begin{center}
\begin{tabular}{@{}ll@{}}
\toprule
Orthonormal plane & Sectional curvature\\
\midrule
$(\xi,E_1)$, $(\xi,E_2)$ & $0$\\
$(\xi,X_j)$, $(\xi,Y_j)$ & $1$\\
$(E_1,E_2)$ & $-nks-2k^2$\\
$(E_a,X_j)$, $(E_a,Y_j)$, $a=1,2$ & $\frac n2(s+2)$\\
$(X_j,Y_j)$ & $-2s-3$\\
$(X_j,X_l)$, $(X_j,Y_l)$, $(Y_j,Y_l)$, $j\ne l$ & $-\frac12s$\\
\bottomrule
\end{tabular}
\end{center}
\medskip
Summing the table we obtain the diagonal components of the Ricci tensor. The off-diagonal components involving only $\xi$, $E_1$, $E_2$ are computed in the same way from the curvature tensor. Namely we have
\begin{gather*}
 \Ric_g(\xi,\xi)=2k,\\
 \Ric_g(E_a,E_a)=-nks-2k^2+kn(s+2)=2k,\\
 \Ric_g(X_j,X_j)=\Ric_g(Y_j,Y_j)
       =1+n(s+2)-2s-3-(k-1)s=2k,\\
 \Ric_g(\xi,E_1)=-\sqrt2kp+2k\frac p{\sqrt2}=0,\\
 \Ric_g(\xi,E_2)=-\sqrt2kq+2k\frac q{\sqrt2}=0,\\
 \Ric_g(E_1,E_2)=2k-2k=0.
\end{gather*}
Finally, all the components $\Ric_g(e_A,e_B)$, $A\neq B$, with at least one of the two vectors in $\KK_0$ vanish identically.
Hence
\begin{equation*}
 \Ric_g=2kg=2(n-1)g.
\end{equation*}

\medskip

Notice that the Reeb vector field is complete. Indeed,
by \eqref{eq:ex-reeb}, its integral curves satisfy
\[
\frac{d}{ds}
\begin{pmatrix}p\\q\end{pmatrix}
=A\begin{pmatrix}p\\q\end{pmatrix},
\qquad
A=\begin{pmatrix}-1&n\\-n&1\end{pmatrix},
\]
together with $t'(s)=1$ and
$x_j'(s)=y_j'(s)=0$.
Since $A^2=-(n^2-1)I$, we have
\[
e^{sA}
=\cos\bigl(\sqrt{n^2-1}\,s\bigr)I
+\frac{\sin\bigl(\sqrt{n^2-1}\,s\bigr)}
{\sqrt{n^2-1}}A.
\]
Thus, for every initial point,
\[
\begin{pmatrix}p(s)\\q(s)\end{pmatrix}
=e^{sA}\begin{pmatrix}p(0)\\q(0)\end{pmatrix},
\qquad t(s)=t(0)+s,
\]
while all $x_j$ and $y_j$ remain constant.
These solutions are defined for every
$s\in\mathbb R$. Hence $\xi$ is complete.
 The metric $g$ is nevertheless incomplete.
Indeed, consider the curve
\[
\gamma(r)=(r,0,\ldots,0),\qquad r\geq0.
\]
As $r\to+\infty$, its first coordinate
tends to infinity, so $\gamma$ does not
converge to any point of $\mathbb R^{2n+1}$.
However, its total Riemannian length is finite, since
\begin{gather*}
 |\dot\gamma(r)|_g^2=
 \frac{n^2+1}{2(n^2-1)^2}
       \left(1+\frac{nr^2}{2(n-1)}\right)^{-2},\\
 L_g(\gamma)=
 \frac{\pi\sqrt{n^2+1}}{2(n+1)\sqrt{n(n-1)}}<\infty.
\end{gather*}
Thus the points $\gamma(r_j)$, for any
sequence $r_j\to+\infty$, form a Cauchy
sequence with respect to the Riemannian
distance, but have no limit in $M$.
Moreover, by the table above, the sectional curvature of the plane spanned by $X_1$ and $Y_1$ equals $-2r^2-3$ at $\gamma(r)$, and hence tends to $-\infty$ along $\gamma$. In particular, the sectional curvature of $g$ is unbounded. 

\medskip

We finally observe that no positive constant
multiple of $g$ is Sasaki--Einstein, even
locally. In fact, suppose by contradiction that
$\widetilde g=cg$ is a
Sasaki--Einstein metric on some nonempty
open subset of $M$,  for some $c>0$. Since
$\Ric_g=2kg$ and the Einstein constant
of a Sasaki--Einstein metric in dimension
$2n+1$ is $2n$, necessarily
\[
c=\frac{k}{n}<1.
\]
Let $\widetilde\xi$ denote the Reeb vector
field of the Sasakian structure associated
with $\widetilde g$. Every Sasakian metric
satisfies
\[
R^{\widetilde g}(A,B)\widetilde\xi
=\widetilde g(B,\widetilde\xi)A
-\widetilde g(A,\widetilde\xi)B.
\]
Since a constant rescaling does not change
the Levi-Civita connection, we have
$R^{\widetilde g}=R^g$ as curvature
tensors of type $(1,3)$. Consequently,
\begin{equation}\label{eq:ex-sasaki-curvature}
R^g(A,B)\widetilde\xi
=c\bigl(g(B,\widetilde\xi)A
-g(A,\widetilde\xi)B\bigr).
\end{equation}
However no nonzero vector
can satisfy this identity. Consider the
Jacobi operator
$\ell=R^g(\,\cdot,\xi)\xi$.
By \eqref{eq:jacobi-tau} and
\eqref{eq:einstein-reeb},
$\ell=P-H+2kJ\tau$.
Since
\[
\tau E_1=E_1,\qquad
\tau E_2=-E_2,\qquad
\tau|_{\KK_0}=0,
\]
we obtain
\begin{equation}\label{eq:ex-jacobi}
\ell E_1=2kE_2,\qquad
\ell E_2=2kE_1,\qquad
\ell|_{\KK_0}=\Id.
\end{equation}
Thus the eigenvalues of $\ell|_{\HH}$
are $2k$, $-2k$, and $1$.
Write
$\widetilde\xi=a\xi+V$,
for $V\in\Gamma(\HH)$.
Taking $A=Z\in\Gamma(\HH)$ and $B=\xi$ in
\eqref{eq:ex-sasaki-curvature}, and
then taking the inner product with
$\xi$, we obtain
$-g(V,\ell Z)=g(R^g(Z,\xi)\widetilde\xi,\xi)=-c\,g(V,Z)$.
Since $\ell$ is self-adjoint, this gives
\[
\ell V=cV.
\]
However, $0<c=k/n<1$, so $c$ is not
an eigenvalue of $\ell|_{\HH}$.
Therefore $V=0$, and hence
$\widetilde\xi=a\xi$, with $a\ne0$.
Now choose a nonzero $Z\in\KK_0$.
Applying \eqref{eq:ex-sasaki-curvature}
once more with $A=Z$ and $B=\xi$,
we obtain
\[
a\ell Z=caZ.
\]
Since $\ell Z=Z$, this implies $c=1$,
contradicting $c=k/n<1$.
Hence no positive constant multiple
of $g$ can be Sasaki--Einstein on
any nonempty open subset of $M$.
This completes the proof of
Theorem~\ref{thm:examples-intro}.


It is worth noticing that the above examples satisfy the divergence identity \eqref{eq:negative-div} explicitly. In the notation of Section~\ref{sec:vanishing}, here $r=1$, $\EE=\operatorname{span}\{E_1,E_2\}$, $\KK=\KK_0$, $SE_1=E_1$ and $SE_2=-E_2$. The structure equations \eqref{eq:ex-structure0}--\eqref{eq:ex-structureK} give
\[
 \LL_{E_1}=\frac1{\sqrt2}(p\Id+qJ)|_{\KK},\qquad
 \LL_{E_2}=\frac1{\sqrt2}(-q\Id+pJ)|_{\KK},
\]
so that the objects defined in \eqref{eq:real-nuW}--\eqref{eq:globalf1} are
\begin{gather*}
 \nu^\sharp=\sqrt2k(pE_1-qE_2),\qquad
 W=\frac{k}{\sqrt2}(pE_1+qE_2),\\
 U=k(p^2+q^2),\qquad f_0=k(p^2-q^2),\qquad f_1=2kpq.
\end{gather*}
Hence the vector field \eqref{eq:negativefield} is
\[
 \mathcal V=\frac{k}{\sqrt2}(qE_1-pE_2)+\frac{k}{2n}(p^2-q^2)\xi.
\]
Since $\divg E_1=\frac{kp}{\sqrt2}$, $\divg E_2=-\frac{kq}{\sqrt2}$ and $\divg\xi=0$, a direct computation based on \eqref{eq:ex-dp}--\eqref{eq:ex-dq} gives
\[
 \divg\mathcal V=-2k-\frac{k}{n}(p^2+q^2),
\]
in accordance with \eqref{eq:negative-div}. Thus it is only the integration over a compact manifold that excludes the intermediate-rank case in Theorem~\ref{thm:main}.

We conclude the section by proving that the only Einstein metric
in the $\mathcal D$-homothetic family is the original metric.

\begin{proposition}
Let $(\phi,\xi,\eta,g)$ be any of the
noncompact Einstein contact metric structures
constructed above. For $a>0$, consider its
$\mathcal D$-homothetic deformation
\[
\eta_a=a\eta,\qquad
\xi_a=a^{-1}\xi,\qquad
\phi_a=\phi,\qquad
g_a=ag+a(a-1)\eta\otimes\eta.
\]
Then $g_a$ is Einstein if and only if $a=1$; more precisely, if $a\neq 1$, then $g_a$ is not Einstein on any nonempty open subset of $\R^{2n+1}$.
Moreover, no positive constant multiple of
$g_a$ is Sasaki--Einstein on any nonempty
open subset.
\end{proposition}

\begin{proof}
Under a $\mathcal D$-homothetic deformation, one has
$\tau_a=a^{-1}\tau$ and $D^a=D$.
The latter identity holds because $D$ is a metric connection also with respect to $g_a$ (recall that $Dg=D\eta=0$) and its torsion is given by \eqref{eq:torsion} also in terms of the deformed structure, since a $g_a$-metric connection is uniquely determined by its torsion.
Since the original metric is Einstein,
Lemma~\ref{lem:einstein} gives
\[
D_\xi\tau=-2nJ\tau.
\]
Consequently,
\[
D^a_{\xi_a}\tau_a
=a^{-2}D_\xi\tau
=-2na^{-2}J\tau.
\]
Suppose that $g_a$ is Einstein on a
nonempty open subset. Since the associated
CR structure remains integrable,
Lemma~\ref{lem:einstein}, applied to the
deformed structure, gives
\[
D^a_{\xi_a}\tau_a
=-2nJ\tau_a=-2na^{-1}J\tau.
\]
Comparing these identities yields
\[
(a-1)J\tau=0.
\]
Since $\tau$ has constant rank two,
we conclude that $a=1$.
Finally, if $a\neq1$, no positive constant multiple of $g_a$ is Einstein on any nonempty open subset, because the Einstein condition is preserved under constant rescalings; in particular, none of them is Sasaki--Einstein. For $a=1$ the claim is part of Theorem~\ref{thm:examples-intro}.
\end{proof}

\section{Final remarks}\label{sec:remarks}

We can derive some immediate consequences of Theorem \ref{thm:main}.

\begin{corollary}
Every compact Einstein contact metric manifold of dimension
$2n+1\ge5$ with integrable associated CR structure has scalar curvature
$2n(2n+1)$. In particular, this class
contains no compact Ricci-flat or negative Einstein metrics.
\end{corollary}

\begin{proof}
Theorem~\ref{thm:main} gives $h=0$, and \eqref{eq:ricxi} then gives
$\lambda=2n$.
\end{proof}

\begin{corollary}
Let $(M^{2n+1},\phi,\xi,\eta,g)$ be a complete Einstein contact metric
manifold with integrable associated CR structure, where $n\ge2$.
If its Einstein constant is positive, then the given structure is
Sasakian.
\end{corollary}

\begin{proof}
The Bonnet--Myers theorem makes $M$ compact. Then we can apply
Theorem~\ref{thm:main}.
\end{proof}

The constancy of the spectrum and the exclusion of maximal rank do
not require compactness. Compactness is used only to integrate
\eqref{eq:negative-div}; the examples of Section~\ref{sec:examples}
show that the conclusion does not follow from the local hypotheses
alone. 

For comparison, in dimension $2n+1\geq 5$, a non-Sasakian contact metric $(\kappa,\mu)$-space (cf. \cite{BKP1995})
satisfies
\[
 h^2=(1-\kappa)P,
\]
with $\kappa<1$, so $h$ has maximal horizontal rank. The Einstein obstruction can also
be seen directly from its Ricci operator \cite[Remark 3.2]{BKP1995}:
\begin{align*}
 Q_g=\left(2(n-1)-n\mu\right)\Id+\left(2(n-1)+\mu\right)h+\left(2(1-n)+n(2\kappa+\mu)\right)\eta\otimes\xi.
\end{align*}
Equality of its eigenvalues on the two nonzero eigenspaces of $h$
forces $\mu=-2(n-1)$ and $\lambda=2(n^2-1)$, whereas the Reeb
direction gives $\lambda=2n\kappa<2n$. This is the same numerical
contradiction as in Proposition~\ref{prop:fullrank}, whose proof
requires no $(\kappa,\mu)$ assumption.

\begin{remark}\label{ex:flat-torus}
The dimension restriction in Theorem~\ref{thm:main} is necessary. In fact, consider the flat torus $T^3=(\R/2\pi\mathbb Z)^3$ with
$g=\dd x^2+\dd y^2+\dd z^2$, and put
\[
 \eta=\cos(2z)\dd x+\sin(2z)\dd y,\qquad
 \xi=\cos(2z)\partial_x+\sin(2z)\partial_y.
\]
Let
\[
 u=-\sin(2z)\partial_x+\cos(2z)\partial_y,\qquad v=\partial_z,
\]
and define $\phi u=v$, $\phi v=-u$, $\phi\xi=0$. Then
\[
 \dd\eta=2\dd z\wedge(-\sin(2z)\dd x+\cos(2z)\dd y)
                              =2g(\,\cdot,\phi\,\cdot).
\]
The $1$-form $\eta$ is a contact form and $\xi$ is its Reeb vector field. The bracket
relations $[\xi,u]=0$ and $[\xi,v]=-2u$ give
$hu=-u$, $hv=v$, so the structure is not Sasakian.
Integrability can also be checked entirely in this real frame.
Its only nonzero Levi--Civita derivatives are
\[
 \nabla_v\xi=2u,\qquad \nabla_vu=-2\xi.
\]
They give $\nabla_\xi\phi=\nabla_u\phi=0$ and
\[
 (\nabla_v\phi)u=0,\qquad
 (\nabla_v\phi)v=2\xi,\qquad
 (\nabla_v\phi)\xi=-2v.
\]
These are exactly \eqref{eq:CRlevi} in view of $hu=-u$ and $hv=v$.
This is a normalization of the classical flat-torus example discussed
in \cite{ADM}. When $n=1$, the full-rank identity
$\lambda=2n^2-2$ gives $\lambda=0$, consistently with this example.
\end{remark}

\end{document}